\documentclass[11pt,letterpaper]{amsart}
\usepackage{mathrsfs,graphicx}
\usepackage{amsopn,amssymb,bbold}

\usepackage{caption}
\usepackage{xy}
\xyoption{all}

\newcommand{\qedlabel}[1]{\renewcommand{\qedsymbol}{\boxnumber{#1}}}
\newcommand{\boxnumber}[1]{{\setlength{\fboxsep}{2pt}\fbox{\textsf{\tiny{#1}}}}}

\DeclareMathOperator{\Hom}{\mathrm{Hom}}
\DeclareMathOperator{\npc}{\mathrm{NPC}}
\DeclareMathOperator{\cat}{\mathrm{CAT}}
\DeclareMathOperator{\hess}{\mathrm{Hess}}
\newcommand{\loc}{\mathrm{loc}}
\newcommand{\DD}{\mathbb{D}}
\newcommand{\half}{\frac{1}{2}}

\newcommand{\B}[1]{\overline{#1}}

\renewcommand{\tilde}{\widetilde}

\newcommand{\tensor}{\otimes}
\newcommand{\into}{\hookrightarrow}
\newcommand{\id}{\mathrm{Id}}
\newcommand{\Z}{\mathbb{Z}}
\newcommand{\R}{\mathbb{R}}
\newcommand{\C}{\mathbb{C}}

\newcommand{\E}{\mathscr{E}}

\newcommand{\boldpoint}[1]{\medskip\par\noindent\textbf{#1}}
\newcommand{\sboldpoint}[1]{\smallskip\par\noindent\textbf{#1}}

\newcommand{\T}{\mathscr{T}}
\renewcommand{\P}{\mathscr{P}}
\newcommand{\AH}{\mathrm{AH}}
\newcommand{\QF}{\mathscr{QF}}
\newcommand{\V}{\mathscr{V}}
\newcommand{\CP}{\mathbb{CP}}
\newcommand{\F}{\mathscr{F}}
\renewcommand{\H}{\mathbb{H}}
\newcommand{\ML}{\mathscr{M\!L}}
\newcommand{\PML}{\mathbb{P} \! \mathscr{M\!L}}
\newcommand{\MF}{\mathscr{M\!F}}
\newcommand{\PQ}{\mathbb{P}^+Q}

\DeclareMathOperator{\PSL}{\mathrm{PSL}}
\DeclareMathOperator{\interior}{Int}
\DeclareMathOperator{\gr}{gr}
\DeclareMathOperator{\pr}{pr}
\DeclareMathOperator{\Gr}{Gr}
\DeclareMathOperator{\area}{Area}
\DeclareMathOperator{\supp}{supp}
\numberwithin{equation}{section}

\theoremstyle{plain}
\newtheorem{thm}{Theorem}[section]
\newtheorem{cor}[thm]{Corollary}
\newtheorem{lem}[thm]{Lemma}

\theoremstyle{remark}
\newtheorem*{rem}{Remark}
\newtheorem*{question}{Question}

\newenvironment{rmenumerate}{\begin{enumerate}
}{\end{enumerate}}

\begin{document}

\title[Schwarzian and Measured Laminations]{The Schwarzian Derivative and Measured Laminations on Riemann Surfaces}
\author{David Dumas}
\date{December 20, 2006}
\address{Department of Mathematics \\ Brown University \\ Providence, RI 02912 \\ USA}
\email{ddumas@math.brown.edu}
\urladdr{http://www.math.brown.edu/\textasciitilde ddumas/}
\thanks{The author was partially supported by an NSF Postdoctoral Research Fellowship}
\begin{abstract}
A holomorphic quadratic differential on a hyperbolic Riemann surface
has an associated measured foliation, which can be straightened
to yield a measured geodesic lamination. On the other hand, a
quadratic differential can be considered as the Schwarzian derivative
of a $\CP^1$ structure, to which one can naturally associate another
measured geodesic lamination using grafting.

We compare these two relationships between quadratic differentials and
measured geodesic laminations, each of which yields a homeomorphism
$\ML(S) \to Q(X)$ for each conformal structure $X$ on a compact
surface $S$.  We show that these maps are nearly the same, differing
by a multiplicative factor of $-2$ and an error term of lower order
than the maps themselves (which we bound explicitly).

As an application we show that the Schwarzian derivative of a $\CP^1$
structure with Fuchsian holonomy is close to a $2\pi$-integral
Jenkins-Strebel differential.  We also study compactifications of the
space of $\CP^1$ structures using the Schwarzian derivative and
grafting coordinates; we show that the natural map between these
extends to the boundary of each fiber over Teichm\"uller space, and we
describe this extension.
\end{abstract}

\maketitle

\section{Introduction}

In this paper we compare two natural homeomorphisms between the vector
space $Q(X)$ of holomorphic quadratic differentials on a compact
Riemann surface $X$ and the $PL$-manifold $\ML(S)$ of measured
laminations on the differentiable surface $S$ underlying $X$.

The first of these two homeomorphisms is a product of $2$-dimensional
geometry--specifically, the theory of measured foliations on Riemann
surfaces.  A holomorphic quadratic differential on $X$ defines a
measured foliation whose leaves are its horizontal trajectories
\cite{Strebel}.  Hubbard and Masur showed that each equivalence class
of measured foliations is obtained uniquely in this way, so that the
resulting map from $Q(X)$ to the space $\MF(S)$ of equivalence classes
of measured foliations is a homeomorphism
\cite{Hubbard-Masur:quadratic-differentials-and-foliations}.

On the other hand, to a measured foliation one can associate a
measured geodesic lamination having the same intersection properties
with simple closed curves.  Thurston showed that the resulting map
$\MF(S) \to \ML(S)$ is a homeomorphism (see
\cite{Levitt:foliations-and-laminations} for a detailed treatment).
Combining this homeomorphism and the Hubbard-Masur map, we obtain a
homeomorphism $\phi_F : \ML(S) \to Q(X)$, which we call the
\emph{foliation map}.

We will compare this map to another homeomorphism between $\ML(S)$ and
$Q(X)$ arising in the theory of complex projective ($\CP^1$)
structures on Riemann surfaces.  Let $P(X)$ denote the space of marked
$\CP^1$ surfaces with underlying Riemann surface $X$.  It is a
classical result that taking the Schwarzian derivatives of the chart
maps induces a homeomorphism $P(X) \to Q(X)$, providing a
complex-analytic parameterization of this moduli space.

Using Thurston's theory of $\CP^1$ structures and a result of Scannell
and Wolf on grafting, one can obtain a more hyperbolic-geometric
description of $P(X)$, leading to a homeomorphism $P(X) \to \ML(S)$.
Here the lamination $\lambda \in \ML(S)$ associated to a projective
structure $Z \in P(X)$ records the bending of a locally convex pleated
plane in $\H^3$ that is the ``convex hull'' of the development of
$\tilde{Z}$ to $\CP^1$ \cite{Kamishima-Tan:grafting}.

As before we combine the two homeomorphisms to obtain a homeomorphism
$\phi_T : \ML(S) \to Q(X)$, which we call the \emph{Thurston map}.

Our goal is to show that in spite of the lack of any apparent
geometric relationship between measured foliations and
complex projective structures, the maps $\phi_T$ and $\phi_F$ are
approximately multiples of one another, up to an error term of smaller
order than either map.  More precisely, we have:

\begin{thm}
\label{thm:schwarzian-nearly-strebel}
Fix a conformal structure $X$ on a compact surface $S$.  Then for all
$\lambda \in \ML(S)$, the foliation map $\phi_F : \ML(S) \to Q(X)$
and the Thurston map $\phi_T : \ML(S) \to Q(X)$ satisfy
\begin{equation*}
\| 2 \phi_T(\lambda) + \phi_F(\lambda) \|_{L^1(X)} 
\leq C(X) \: \left ( 1 + \sqrt{\|\phi_F(\lambda) \|_{L^1(X)}} \right )
\end{equation*}
where $C(X)$ is a constant that depends only on $X$.
\end{thm}

\boldpoint{The Schwarzian derivative.}  The constructions used to
define the foliation map $\phi_F$ and the Thurston parameterization of
$\CP^1$ structures are essentially geometric, and can be understood as
part of the rich interplay between the theory of Riemann surfaces and
hyperbolic geometry in two and three dimensions.  The Schwarzian
derivative seems more analytic than geometric, however, and relating the
quadratic differential obtained from the charts of a $\CP^1$ structure
(using the Schwarzian) to geometric data accounts for a significant
part of the proof of Theorem \ref{thm:schwarzian-nearly-strebel}.

Relationships between the Schwarzian derivative and geometry have been
explored by a number of authors, using curvature of curves in the
plane \cite{Flanders:schwarzian-as-curvature}
\cite{CDO:curvature-planar-harmonic}, curvature of surfaces in
hyperbolic space \cite{Epstein:hyperbolic-gauss-map}, Lorentzian
geometry \cite{DO:world-lines-and-schwarzian} \cite{OT:sturm-theory},
and osculating M\"obius transformations \cite{Thurston:zippers}
\cite{An}.  Furthermore, there are a number of different ways to
generalize the classical Schwarzian derivative to other kinds of maps,
more general domain and range spaces, or both (e.g.
\cite{Ahlfors:cross-ratios-and-schwarzian}
\cite{GF:schwarzian-in-several-complex-variables}
\cite{BO:cocycles-diff-s1}).

The proof of Theorem \ref{thm:schwarzian-nearly-strebel} turns on the
decomposition of the Schwarzian derivative of a $\CP^1$ structure into
a sum of two parts (Theorem \ref{thm:schwarzian-decomposition}), one
of which is manifestly geometric, and another which is ``small'',
having norm bounded by a constant depending only on $X$.  This
decomposition is a product of the Osgood-Stowe theory of the
Schwarzian \emph{tensor}--a particular generalization of the
Schwarzian derivative that measures the difference between two
conformally equivalent Riemannian metrics \cite{OS:schwarzian}.  Much
as the Schwarzian derivative of a composition of holomorphic maps can
be decomposed using the cocycle property (or ``chain rule''), the
Schwarzian tensor allows us to decompose the quadratic differential
associated to a $\CP^1$ structure by finding a conformal metric that
interpolates those of the domain and range of the projective charts.
In our case, the appropriate interpolating metric is the
\emph{Thurston metric} associated to a grafted surface, which is a
sort of Kobayashi metric in the category of $\CP^1$ surfaces (see
\S\ref{sec:grafting-differential}).

We apply this decomposition of the Schwarzian to $\phi_T(\lambda)$,
the Schwarzian derivative of the $\CP^1$ structure on $X$ with
grafting lamination $\lambda \in \ML(S)$.  The two terms from this
decomposition of are then analyzed separately, ultimately leading to a
proof of Theorem \ref{thm:schwarzian-nearly-strebel}.

\boldpoint{Harmonic maps.}  The first part of the decomposition of
$\phi_T(\lambda)$ is the \emph{grafting differential} $\Phi(\lambda)$,
a non-holomorphic quadratic differential on $X$ whose trajectories
foliate the grafted part of the $\CP^1$ structure.  The grafting
differential can also be defined as (a multiple of) the Hopf
differential of the \emph{collapsing map} of the $\CP^1$ structure,
which is a generalization of the retraction of a set in $\CP^1$ to the
boundary of the hyperbolic convex hull of its complement.

Using harmonic maps techniques, one can bound the difference between
the grafting differential and Hubbard-Masur differential
$\phi_F(\lambda)$ in terms of the extremal length of $\lambda$ (see
\S\ref{sec:grafting-differential}, also \cite{Dumas:erratum}
\cite{Dumas:antipodal}):
\begin{equation*}
\| \Phi(\lambda) - \phi_F(\lambda) \|_{L^1(X)} 
\leq C\:(1 + E(\lambda,X)^\frac{1}{2})
\end{equation*}
Here the constant $C$ depends only on the topological type of $X$ (see
\S\ref{sec:grafting-differential}).  Geometrically, this makes sense:
the grafting differential is holomorphic where it is nonzero, which is
most of $X$ when $\lambda$ is large (see Figure
\ref{fig:thurston-metric} in \S\ref{sec:thurston-metric}); thus it
should be close to the unique holomorphic differential
$\phi_F(\lambda)$ with the same trajectory structure.

\boldpoint{Analytic methods.}  The second part of the decomposition of 
$\phi_T(\lambda)$ is the Schwarzian tensor of the Thurston metric of
the $\CP^1$ structure relative to the hyperbolic metric of $X$.  We
complete the proof of Theorem \ref{thm:schwarzian-nearly-strebel} by
showing that this non-holomorphic quadratic differential is bounded,
using analysis and the properties of the Thurston metric.

Essential to this bound is the fact that the Schwarzian tensor of a
conformal metric (relative to any fixed background metric) depends
only on the $2$-jet of its density function, and thus the norm of this
tensor is controlled by the norm of the density function in an
appropriate Sobolev space.  Using standard elliptic theory, it is
actually enough to bound the Laplacian of the density function, which
is essentially the curvature $2$-form of the Thurston metric.

On the other hand, the curvature of the Thurston metric can be
understood using the measured lamination $\lambda$ and the properties
of grafting.  We ultimately show that the Thurston metric's curvature
is concentrated near a finite set of points, except for an
exponentially small portion that may diffuse into the rest of $X$.
This is enough to bound the Schwarzian tensor on a hyperbolic disk of
definite size, and then, using a compactness argument, on all of $X$.

\boldpoint{Applications.}  As an application of Theorem
\ref{thm:schwarzian-nearly-strebel}, we show that the countably many
$\CP^1$ structures on $X \in \T(S)$ with Fuchsian holonomy are close
to the $2 \pi$-integral Jenkins-Strebel differentials in $Q(X)$
(Theorem \ref{thm:fuchsian-centers}).  In particular these \emph{Fuchsian
centers} are arranged in a regular pattern, up to an error term of
smaller order than their norms (though the difference may be unbounded
for a sequence of centers going to infinity).  Numerical experiments
illustrating this effect are presented in
\S\ref{sec:applications-fuchsian}.

We also apply the main theorem to the fiberwise compactification of
the space of $\CP^1$ structures induced by the Schwarzian map $P(X)
\xrightarrow{\sim} Q(X)$, and show that the map to the grafting
coordinates $P(X) \into \ML(S) \times \T(S)$ has a natural
continuous extension (Theorem \ref{thm:schwarzian-compactification}).
This extends the results of \cite{Dumas:antipodal}, where this
boundary of $P(X)$ was studied independently of the Schwarzian
parameterization using the \emph{antipodal involution} $i_X : \PML(S)
\to \PML(S)$.

\boldpoint{Outline of the paper.}
\sboldpoint{\S\ref{sec:conformal-metrics}} presents some
background material on conformal metrics and tensors on Riemann
surfaces, which are the main objects of study in the proof of Theorem
\ref{thm:schwarzian-nearly-strebel}.

\sboldpoint{\S\ref{sec:hubbard-masur}} describes the Hubbard-Masur
construction of a homeomorphism between $Q(X)$ and $\MF(S)$, and the
associated \emph{foliation map} $\phi_F : \ML(S) \to Q(X)$.

\sboldpoint{\S\ref{sec:grafting-and-cp1-structures}} describes
another connection between $\ML(S)$ and $Q(X)$ using Thurston's theory
of grafting and $\CP^1$ structures on Riemann surfaces and results of
Scannell-Wolf on grafting.  The result is the \emph{Thurston map}
$\phi_T : \ML(S) \to Q(X)$.

\sboldpoint{\S\ref{sec:grafting-differential}} continues our
study of grafting and projective structures by introducing the
\emph{Thurston metric} and the \emph{grafting differential}.  These
natural geometric objects are analogous to the singular Euclidean
metric and the holomorphic quadratic differential that arise in the
Hubbard-Masur construction.

\sboldpoint{\S\ref{sec:schwarzian-derivative}} introduces the
Osgood-Stowe Schwarzian tensor, a generalization of the Schwarzian
derivative that we use to relate the Schwarzian of a projective
structure to the analytic properties of conformal metrics.

\sboldpoint{\S\ref{sec:decomposition-of-the-schwarzian}} establishes a
  decomposition of the Schwarzian of a projective structure (Theorem
  \ref{thm:schwarzian-decomposition}) into the sum of the grafting
  differential and the Schwarzian tensor of the Thurston metric.  This
  decomposition is the main conceptual step toward the proof of the
  main theorem.

\sboldpoint{\S\S\ref{sec:npc-metrics}-\ref{sec:regularity-and-compactness}}
  begin our analytic study of the Thurston metric by deriving certain
  regularity properties of nonpositively curved conformal metrics on
  Riemann surfaces (of which the Thurston metric is an example).

\sboldpoint{\S\S\ref{sec:thurston-metric}-\ref{sec:schwarzian-of-thurston-metric}}
contain the key estimates on the Thurston metric, showing that its
curvature concentrates near a finite set of points (Theorem
\ref{thm:curvature-concentration}) away from which its Schwarzian
tensor is bounded (Theorem
\ref{thm:schwarzian-of-thurston-metric-bounded}).

\sboldpoint{\S\ref{sec:proof}} contains the proof Theorem
\ref{thm:schwarzian-nearly-strebel}, which combines the bound on the
Schwarzian tensor of the Thurston metric, the decomposition of the
Schwarzian of a projective structure, and properties of the grafting
differential (Theorems
\ref{thm:schwarzian-of-thurston-metric-bounded},
\ref{thm:schwarzian-decomposition}, and
\ref{thm:grafting-differentials-close}, respectively).

\sboldpoint{\S\S\ref{sec:applications-fuchsian}-\ref{sec:applications-compactification}}
present the applications of Theorem
\ref{thm:schwarzian-nearly-strebel} mentioned above---locating Fuchsian
centers in $P(X)$ (and related numerical experiments), and the
continuous extension of the grafting coordinates $P(X) \into \ML(S)
\times \T(S)$ to respective compactifications.

\boldpoint{Acknowledgements.}  The author thanks Georgios
Daskalopoulos, Bob Hardt, Curt McMullen, and Mike Wolf for stimulating
discussions related to this work.  He is also grateful to the
anonymous referee for several suggestions that improved the paper.
Some of this work was completed while the author was a postdoctoral
fellow at Rice University, and he thanks the the department for its
hospitality.

\boldpoint{Notational conventions.}  In what follows, $X$ denotes a
compact (except in \S\ref{sec:applications-fuchsian}) hyperbolic
Riemann surface and $S$ the underlying differentiable surface of genus
$g$ and Euler characteristic $\chi = 2 - 2g$.

The expression $C(a,b,\ldots)$ is used to indicate that an unspecified
constant $C$ depends on quantities $a, b, \ldots$, one of which is
typically a conformal structure $X$.

\section{Conformal metrics and tensors}
\label{sec:conformal-metrics}

We briefly recall some constructions related to conformal metrics and
tensors on Riemann surfaces that will be used in the sequel.  On a
fixed compact Riemann surface $X$, choose a complex line bundle
$\Sigma$ with $\Sigma^2 = K_X$.  This allows us to define the bundle
$S_{i,j} = S_{i,j}(X)$ of differentials of type $(i,j)$ for all $(i,j)
\in (\half \Z)^2$:
\begin{equation*}
S_{i,j} = \Sigma^{2 i} \bar{\Sigma}^{2 j}
\end{equation*}
We will be most interested in $S_{\half,\half}$, whose sections
include \emph{conformal metrics} on $X$, which in a coordinate chart
have the form
\begin{equation*}
\rho(z) \: |dz|
\end{equation*}
where $\rho$ is a nonnegative function.  Such a conformal metric
determines length and area functions,
\begin{equation*}
\begin{split}
\ell(\gamma,\rho) = \int_\gamma \rho & \;\; 
\text{ where } \; \gamma : [0,1] \to X\\
A(\Omega,\rho) = \int_\Omega \rho^2 & \;\; 
\text{ where } \; \Omega \subset X
\end{split}
\end{equation*}
which, modulo sufficient regularity and positivity of $\rho$, make $X$
into a geodesic metric space.  The properties of a particular class of
these metrics will be studied further in
\S\S\ref{sec:npc-metrics}-\ref{sec:regularity-and-compactness}.

The $L^p$ norm of a section of $S_{i,j}$ is well-defined independent
of any background metric on $X$ when $p(i+j) = 2$.  For example, the
$L^2$ norm on $S_{\half,\half}$ corresponds to the area of a conformal
metric.  Given a conformal metric $\rho \in L^2(S_{\half,\half})$, we
can also define the $L^p$ norm of a section $\xi$ of $S_{i,j}$ with
respect to $\rho$:
$$ \| \xi \|_{L^p(S_{i,j}, \rho)} = \int_X |\xi|^p \rho^{2 - p(i+j)}.$$
When the tensor type is understood and a background metric $\rho$ is
fixed (or unnecessary) we will abbreviate this norm as $\| \xi
\|_{L^p(X)}$ or simply $\|\xi\|_p$.

A pullback construction will provide most of our examples of conformal
metrics; specifically, given a smooth map $f : X \to (M,g)$, where $M$
is a Riemannian manifold, the pullback metric $f^*(g)$ need not lie in
the conformal class of $X$, however it can be decomposed relative to
this conformal structure as follows:
\begin{equation*}
\begin{split}
f^*(g) &= \Phi(f) + \rho(f)^2 + \B{\Phi(f)}\\
\end{split}
\end{equation*}
Here $\Phi(f) \in \Gamma(S_{2,0})$ is the \emph{Hopf differential} of
$f$ and $\rho(f) \in \Gamma(S_{\half,\half})$ is a conformal metric
that is the ``isotropic part'' of the pullback of the line element of
$g$.  Of course the smoothness assumption on $f$ may be relaxed
considerably; the natural regularity class for our purposes is the
space $W^{1,2}(X,M)$ of Sobolev maps into $M$ with $L^2$ derivatives.
For such maps, the resulting conformal metric $\rho$ is in
$L^2(S_{\half,\half})$ and the Hopf differential $\Phi$ lies in
$L^1(S_{2,0})$.

\section{The Hubbard-Masur construction}
\label{sec:hubbard-masur}

In this section we review the first of two relationships between
measured laminations and quadratic differentials that we will
explore.  A result of Hubbard and Masur is essential to
this relationship.

Let $Q(X) \subset L^1(S_{2,0}(X))$ denote the space of holomorphic
quadratic differentials on $X$.  Each differential $\phi \in Q(X)$ has
an associated (singular) measured foliation $\F(\phi)$ whose leaves
integrate the distribution of vectors $v \in TX$ satisfying $\phi(v)
\geq 0$.  There is also a well-known homeomorphism $\MF(S) \simeq
\ML(S)$ between the spaces of (measure equivalence classes of)
measured foliations and measured geodesic laminations on a hyperbolic
surface.  Roughly speaking, to obtain a lamination from a measured
foliation, one replaces the nonsingular leaves of the foliation with
geodesic representatives for the hyperbolic metric on $X$ (for a
detailed account, see \cite{Levitt:foliations-and-laminations}).  Thus
we obtain a map $\Lambda : Q(X) \to \ML(S)$.

Hubbard and Masur showed that $\F$ (and thus also $\Lambda$) is a
homeomorphism
\cite{Hubbard-Masur:quadratic-differentials-and-foliations}; in other
words, given any measure equivalence class of measured foliations
$\F_0$, there is a unique holomorphic quadratic differential $\phi$
with $\F(\phi) \sim \F_0$.  (For other perspectives on this result see
\cite{Kerckhoff:asymptotic-geometry}
\cite{Gardiner:measured-foliations}
\cite{Wolf:realizing-measured-foliations}, and for the special case of
foliations with closed trajectories see
e.g. \cite{Jenkins:extremal-metrics}
\cite{Strebel:closed-trajectories} \cite{Gardiner:Jenkins-Strebel}
\cite{Marden-Strebel:heights} \cite{Wolf:harmonic-maps-to-graphs}.)
We will be interested in the inverse homeomorphism
\begin{equation*}
\phi_F : \ML(S) \to Q(X)
\end{equation*}
which we call the \emph{foliation map}, as it associates to $\lambda
\in \ML(S)$ a quadratic differential whose foliation has prescribed
measure properties.

The foliation map $\phi_F : \ML(S) \to Q(X)$ is well-behaved with
respect to natural structures on the spaces $\ML(S)$ and $Q(X)$; it
preserves basepoints (i.e. $\phi_F(0) = 0$, where $0 \in \ML(S)$ is
the empty lamination), and is homogeneous of degree $2$ on rays in
$\ML(S)$.

\section{Grafting and $\CP^1$ structures}
\label{sec:grafting-and-cp1-structures}

One might say that the foliation homeomorphism $\phi_F : \ML(S) \to
Q(X)$ exists because both $\ML(S)$ and $Q(X)$ are models for the space
$\MF(S)$ of measured foliation classes on $S$; one model comes from
hyperbolic geometry ($\ML(S) \simeq \MF(S)$), while the other comes
from the singular Euclidean geometry of a quadratic differential
($Q(X) \simeq \MF(S)$).

In a similar vein, we now construct a homeomorphism $\phi_T : \ML(S) \to
Q(X)$, which we call the \emph{Thurston map}, by showing that both
$\ML(S)$ and $Q(X)$ are naturally homeomorphic to the space $P(X)$ of
complex projective structures on $X$.  The identification of $\ML(S)$
with $P(X)$ will involve hyperbolic geometry (in $\H^3$) and
Thurston's projective version of \emph{grafting}, while that of $Q(X)$
with $P(X)$ uses the Schwarzian derivative.  We begin with a few
generalities on $\CP^1$-structures.

A complex projective structure on $S$ is an atlas of charts with
values in $\CP^1$ and M\"obius transition functions.  The space $\P(S)$
of marked $\CP^1$ structures fibers over Teichm\"uller space by the map
$\pi : \P(S) \to \T(S)$ which records the underlying complex structure.
When the underlying complex structure of a $\CP^1$ surface is $X$, we
say it is a $\CP^1$ structure \emph{on $X$}.

Let $Z$ be a projective structure on $X$.  The chart maps of $Z$ can
be analytically continued on the universal cover $\tilde{X} \simeq
\DD$ to give a locally univalent holomorphic map $f : \tilde{X}
\to \CP^1$, called the \emph{developing map}.  This map is not unique,
but any two such maps differ by composition with a M\"obius
transformation.  When restricted to any open set on which it is
univalent, the developing map is a projective chart.

Because the projective structure on $\tilde{X}$ induced by lifting
$Z$ is invariant under the action of $\pi_1(S)$ by deck
transformations, for any $\gamma \in \pi_1(S)$ and $z \in
\tilde{X}$, the germs $f_{z}$ and $f_{\gamma z}$ differ by
composition with a M\"obius transformation $A_\gamma$.  The map
$\gamma \mapsto A_\gamma$ defines a homomorphism $\eta(Z) : \pi_1(S)
\to \PSL_2(\C)$, called the \emph{holonomy representation} of $Z$,
which is unique up to conjugation.

To obtain a concrete realization of the fiber $P(X) = \pi^{-1}(X)$, we
use the \emph{Schwarzian derivative}, a M\"obius-invariant differential
operator on locally injective holomorphic maps to $\CP^1$:
\begin{equation*}
S(f) = \left [ \left(\frac{f''(z)}{f'(z)}\right)' -
\half \left(\frac{f''(z)}{f'(z)}\right)^2 \right ] \: dz^2
\end{equation*}
By its M\"obius invariance, the Schwarzian derivatives of the charts
of a $\CP^1$ structure on $X$ (viewed as maps from subsets of $\DD
\simeq \H^2 \simeq \tilde{X}$ to $\CP^1$) join together to form a
holomorphic quadratic differential on $X$.  Equivalently, the
Schwarzian derivative of the developing map $S(f) \in Q(\DD)$ is
invariant under the action of $\pi_1(S)$ by deck transformations of
the universal covering $\DD \to X$, and so it descends to a quadratic
differential on $X$.  The resulting map
\begin{equation*}
\xymatrix{
P(X) \ar[r]^-{S} & Q(X)
}
\end{equation*}
is a homeomorphism, and in fact a biholomorphism with respect to the
natural complex structure of $P(X)$.  This is the \emph{Poincar\'e
parameterization} of $P(X)$ (see \cite{Hejhal}, \cite{Gunning}).

\begin{figure}
\centerline{\includegraphics{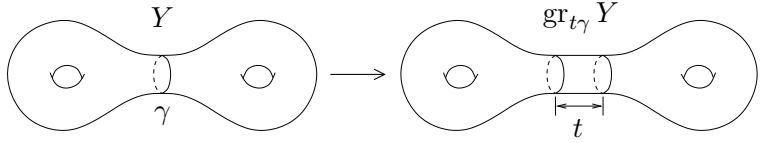}}
\caption{The basic example of grafting.}
\label{fig:basicgraft}
\end{figure}

An alternate and more geometric approach to $\CP^1$ structures was
described by Thurston using \emph{grafting}, a cut-and-paste operation
on hyperbolic Riemann surfaces (see e.g. \cite{Ma69} \cite{He75}
\cite{Kamishima-Tan:grafting} \cite{SW:grafting}). The conformal
grafting map
\begin{equation*}
\gr : \ML(S) \times \T(S) \to \T(S)
\end{equation*} sends the pair $(\lambda,Y)$ to a
surface obtained by removing the geodesic lamination supporting
$\lambda$ from $Y$ and replacing it with a ``thickened lamination''
that has a Euclidean structure realizing the measure of $\lambda$.
The details of the construction are more easily explained when
$\lambda$ is supported on a simple closed geodesic $\gamma$ with
weight $t$, in which case $\gamma$ is simply replaced with the
Euclidean cylinder $\gamma \times [0,t]$ to obtain the grafted surface
$\gr_{t \gamma} Y$ (see Figure \ref{fig:basicgraft}).  As such
weighted curves are dense in $\ML(S)$, the existence of the grafting
operation in general can be reduced to a continuity property; for
details, see \cite{Kamishima-Tan:grafting}.

Thurston introduced a projective grafting map $\Gr : \ML(S) \times
\T(S) \to \P(S)$ that puts a canonical projective structure on a
grafted surface.  This map associates to $(\lambda,Y)$ a projective
structure on $\gr_\lambda Y$ whose ``convex hull boundary'' in $\H^3$
is the locally convex pleated surface obtained by bending $\tilde{Y}
\simeq \H^2$ along the lift of the measured lamination $\lambda$.
Roughly speaking, the Gauss map from this surface (which follows
normal geodesic rays out to $\CP^1 = \partial_\infty\H^3$) provides a
system of charts for this projective structure.  This projective
version of grafting is especially interesting because every projective
structure can be obtained from grafting in exactly one way; that is,
\begin{equation*}
\xymatrix{
\ML(S) \times \T(S) \ar[r]^-{\Gr} & \P(S)
}
\end{equation*}
is a homeomorphism.  This is \emph{Thurston's theorem}, a detailed
proof of which can be found in \cite{Kamishima-Tan:grafting}.

Thus grafting gives a geometric description of $\P(S)$, and we can
associate to any projective surface a pair $(\lambda,Y) \in \ML(S)
\times \T(S)$.  However, it is not immediately clear how the fibers
$P(X)$ with a fixed underlying conformal structure fit into this
grafting picture.  We will now describe how one of the grafting
coordinates, the measured lamination $\lambda \in \ML(S)$, suffices to
parameterize any such fiber.

Scannell and Wolf showed that for each $\lambda \in \ML(S)$ the
conformal grafting map $\gr_\lambda : \T(S) \to \T(S)$ is a
homeomorphism \cite{SW:grafting}; we call the inverse homeomorphism
\emph{pruning} by $\lambda$:
\begin{equation}
\label{eqn:pruning}
\pr_\lambda = \gr_\lambda^{-1} : \T(S) \to \T(S).
\end{equation}
Thus for each lamination $\lambda \in \ML(S)$, the hyperbolic surface
\begin{equation*}
Y_\lambda = \pr_\lambda X
\end{equation*}
has the property that when grafted by $\lambda$, the resulting surface
is conformally equivalent to $X$.  If we use Thurston's projective
extension of grafting, then the result of grafting $Y_\lambda$ by
$\lambda$ is a projective structure on $X$, which we call
$X(\lambda)$, i.e.
\begin{equation*}
X(\lambda) = \Gr_\lambda \pr_\lambda X \in P(X).
\end{equation*}
The resulting map
\begin{equation*}
\begin{split}
\beta: \ML(S) &\to P(X)\\
\lambda &\mapsto X(\lambda)
\end{split}
\end{equation*}
is evidently a homeomorphism, because its inverse is the composition
of $\Gr^{-1}|_{P(X)} : P(X) \to \ML(S) \times \T(S)$ with the
projection of $\ML(S) \times \T(S)$ onto the first factor.

By taking the Schwarzian of the projective structure $X(\lambda)$, we
obtain the \emph{Thurston map} $\phi_T : \ML(S) \xrightarrow{\sim}
Q(X)$:
\begin{equation*}
\xymatrix@C+5mm{
\ML(S) \ar@/_1pc/[rr]_-{\phi_T} \ar[r]^-{\beta} 
& P(X) \ar[r]^-{S} & Q(X)
}
\end{equation*}

When compared to the foliation map, the Thurston map $\phi_T$ is
somewhat opaque; it is a homeomorphism, and it is easily seen to map
the empty lamination to the zero differential (i.e. $\phi_T(0)=0$),
but there is no reason expect this map to be homogeneous, or even to
preserve rays (see \S\ref{sec:applications-fuchsian} for numerical
experiments suggesting that it does not).  The main obstruction to an
intuitive understanding of this map would seem to be the lack of a
connection between the analytic definition of the Schwarzian
derivative and geometric properties of the projective surface.  To put
it another way, the Schwarzian $\phi$ of a projective structure on $X$
has an associated foliation $\F(\phi)$, but there is no obvious
relationship between the geometry of this foliation and that of the
$\CP^1$ structure.

It is just such a relationship we hope to reveal by relating the
Hubbard-Masur construction (and $\phi_F$) to the Thurston map
$\phi_T$, if only approximately.

\section{The Thurston metric and grafting differential}
\label{sec:grafting-differential}

The foliation map $\phi_F : \ML(S) \to Q(X)$ associates a holomorphic
quadratic differential $\phi_F(\lambda)$ to each measured lamination
$\lambda$, which in turn gives a singular Euclidean conformal metric
$|\phi_F(\lambda)|^{1/2}$ on $X$.  This metric is the natural one in
which to examine the measured foliation $\F(\phi_F(\lambda))$, and is
extremal for this foliation class in the sense of extremal length.

We now start to develop a similar picture for the Thurston map.  In
this case each measured lamination $\lambda$ gives rise to a
projective structure $X(\lambda) \in P(X)$ of the form $\Gr_\lambda
Y_\lambda$ whose Schwarzian is $\phi_T(\lambda)$.  The grafting
construction presents $X(\lambda)$ as a union of two parts: The
hyperbolic part, $X_{-1}(\lambda)$, which comes from $Y_\lambda$, and
a Euclidean part, $X_{0}(\lambda)$, which is grafted into $Y_\lambda$
along the geodesic lamination $\lambda$.

There is a natural conformal metric $\rho_\lambda \in
L^2(S_{\half,\half}(X))$ on $X$ associated to the projective structure
$X(\lambda)$, called the \emph{Thurston metric}, which combines
the hyperbolic metric on $X_{-1}(\lambda)$ and the Euclidean metric on
$X_0(\lambda)$ (see \cite[\S2.2]{SW:grafting},
\cite{KP:canonical-metric}).  For example, $\rho_0$ is the hyperbolic
metric on $X$.

Another convenient description of the Thurston metric is as the
``Kobayashi metric'' in the category of $\CP^1$ surfaces and locally
M\"obius maps \cite[\S2.1]{Tanigawa:harmonic-maps}.  To make this
precise, we define a projectively immersed disk in a $\CP^1$ surface
$Z$ to be a locally M\"obius map $\delta$ from the unit disk in $\C$
(with its canonical $\CP^1$ structure) to $Z$.  Then the Thurston
length of a tangent vector $v \in T_x Z$ is the minimum length it is
assigned by the hyperbolic metric of $\DD$ when pulled back via a
projective immersion $\delta : \DD \to Z$ with $\delta(0) = x$.

Since the hyperbolic metric $\rho_0$ is the ordinary Kobayashi metric
for $X$, where length is obtained as an infimum over the class of
holomorphic immersions of the disk (which is larger than the class of
projective immersions), we conclude immediately from this definition
that for all $\lambda \in \ML(S)$,
\begin{equation}
\label{eqn:thurston-larger-than-hyperbolic}
 \rho_\lambda \geq \rho_0.
\end{equation}

The Thurston metric is also related to the \emph{collapsing map}
$\kappa : X \to Y_\lambda$, which collapses the grafted part of a
surface orthogonally onto the geodesic representative of $\lambda$ on
$Y_\lambda$.  This map is distance non-increasing for the Thurston
metric on $X$, and the Thurston metric is pointwise the smallest
conformal metric on $X$ with this property (since at every point there
is a direction in which $\kappa$ is an isometry with respect to the
Thurston metric \cite[Thm.~8.6]{KP:canonical-metric}, any smaller
conformal metric on $X$ would make the map expand somewhere).

The Hopf differential $\Phi(\kappa) \in L^2(S_{2,0}(X))$ of the
collapsing map has an associated partial measured foliation
$\F(\Phi(\kappa))$ that is supported in the grafting locus
$X_0(\lambda)$, and whose leaves are Euclidean geodesics.  Let us
define the \emph{grafting differential} $\Phi(\lambda)$,
\begin{equation}
\label{eqn:grafting-differential}
\Phi(\lambda) = 4 \Phi(\kappa).
\end{equation}
The normalization is chosen so that the partial measured foliation
$\F(\Phi(\lambda))$ represents the measure equivalence class $\lambda$
(see \cite[\S6]{Dumas:antipodal}), i.e.
\begin{equation*}
\F(\Phi(\lambda)) \sim \lambda,
\end{equation*}
where $a \sim b$ means that $a$ and $b$ have the same intersection
numbers with all simple closed curves.  Thus $\Phi(\lambda)$ and
$\F(\Phi(\lambda))$ are to the Thurston metric much as a holomorphic
differential $\phi$ and foliation $\F(\phi)$ are to the singular
Euclidean metric $|\phi|^{1/2}$.

Much about the large-scale behavior of grafting and related objects
can be determined using the fact that the collapsing map $\kappa : X
\to Y_\lambda$ is nearly harmonic, i.e. it nearly minimizes energy in
its homotopy class.

\begin{thm}[{Tanigawa, \cite[Thm.~3.4]{Tanigawa:harmonic-maps}}]
\label{thm:grafting-energy-close}
Let $X = \gr_\lambda Y$ with collapsing map $\kappa : X \to Y$, and
let $h : X \to Y$ denote the harmonic map compatible with the markings
of $X$ and $Y$, and $\E(h)$ its energy.  Then
\begin{equation*}
\E(h) \approx \E(\kappa) \approx \half \ell(\lambda,Y) \approx
\half E(\lambda,X) = \half \| \phi_F(\lambda) \|_1
\end{equation*}
where $E(\lambda,X)$ is the extremal length of $\lambda$ on $X$ and
$\ell(\lambda,Y)$ is the length of $\lambda$ with respect to the
hyperbolic metric on $Y$.  Here $A \approx B$ means that the
difference $A - B$ is bounded by a constant that depends only on the
topology of $X$.
\end{thm}

\begin{rem}
In \cite{Tanigawa:harmonic-maps}, Tanigawa establishes a set of
inequalities relating $\E(\kappa) = \frac{1}{2} \ell(\lambda,Y) + 2
\pi |\chi(S)|$, $\E(h)$, and $E(\lambda,X)$, from which the
approximate equalities in Theorem \ref{thm:grafting-energy-close}
follow by algebra.  See \cite[\S7]{Dumas:antipodal} for details.
\end{rem}

Harmonic maps techniques can also be used to relate the grafting
differential to the foliation map; in fact, we have:

\begin{thm}[{\cite[Thm. 10.1]{Dumas:erratum},  \cite{Dumas:antipodal}}]
\label{thm:grafting-differentials-close}
For any $X \in \T(S)$ and $\lambda \in \ML(S)$, the holomorphic
quadratic differential $\phi_F(\lambda) \in Q(X)$ and the grafting
differential $\Phi(\lambda) \in L^1(S_{2,0}(X))$ satisfy
\begin{equation*}
\| \Phi(\lambda) - \phi_F(\lambda) \|_1 \leq 
C \left ( 1 + E(\lambda,X)^\half \right ).
\end{equation*}
\end{thm}

Since Theorem \ref{thm:grafting-differentials-close} has an important
role in the proof of Theorem \ref{thm:schwarzian-nearly-strebel}, we
take a moment to sketch the ideas behind it: First, a construction
dual to that of the collapsing map gives a \emph{co-collapsing map}
$\hat{\kappa} : \Tilde{X} \to T_\lambda$ with Hopf differential
$\Phi(\hat{\kappa}) = -\frac{1}{4} \Phi(\lambda)$, where $T_\lambda$
is the $\R$-tree dual to $\lambda$.  An energy estimate in the spirit
of Theorem \ref{thm:grafting-energy-close} shows that the
co-collapsing map is nearly harmonic, and by a theorem of Wolf, the
Hopf differential of the harmonic map to $T_\lambda$ is
$-\frac{1}{4}\phi_F(\lambda)$.  Finally, an estimate of
Korevaar-Schoen from \cite[\S 2.6]{KS:sobolev-spaces} shows that a
nearly-harmonic map to a tree (or indeed, any $CAT(0)$ metric space)
has Hopf differential which is close to that of the harmonic map,
leading to the specific bound in Theorem
\ref{thm:grafting-differentials-close}.

What is missing from this harmonic maps picture is any geometric
control on the Schwarzian $\phi_T(\lambda)$.  In the next section we
discuss the Osgood-Stowe generalization of the Schwarzian derivative,
which we then use in \S\ref{sec:decomposition-of-the-schwarzian} to
relate the Schwarzian and the grafting differential.

\section{The Schwarzian derivative and Schwarzian tensor}
\label{sec:schwarzian-derivative}

In a precise sense, the Schwarzian derivative measures the extent to
which a locally injective holomorphic map $f$ fails to be a M\"obius
transformation \cite{Thurston:zippers}; for example, $f$ is (the
restriction of) a M\"obius transformation if and only if $S(f) = 0$.
The Schwarzian of a composition of maps is governed by the
\emph{cocycle relation}:
\begin{equation}
\label{eqn:schwarzian-cocycle}
S(f \circ g) = g^* S(f) + S(g)
\end{equation}

In \cite{OS:schwarzian}, Osgood and Stowe construct a generalization of
the Schwarzian derivative that acts on a pair of conformally
equivalent Riemannian metrics on a manifold.  We describe this
generalization only in the case of conformal metrics on a Riemann
surface, as this is the case we will use.

Given conformal metrics $\rho_1,\rho_2 \in L^2(S_{\half,\half}(X))$,
define
\begin{equation}
\label{eqn:sigma}
\sigma(\rho_1,\rho_2) = \log(\rho_2 / \rho_1).
\end{equation}
Then the \emph{Schwarzian tensor $\beta(\rho_1,\rho_2)$ of $\rho_2$
relative to $\rho_1$} is defined as
\begin{equation}
\label{eqn:hessian-definition}
\beta(\rho_1,\rho_2) = [\hess_{\rho_1} (\sigma) - d\sigma \tensor
 d\sigma]^{2,0}
\end{equation}
where we have written $\sigma$ instead of $\sigma(\rho_1,\rho_2)$ for
brevity.  This definition differs from that of Osgood and Stowe in
that we take only the $(2,0)$ part, whereas they consider the
traceless part, which in this case is the sum of $\beta$ and its
complex conjugate\footnote{We also use a slightly different notation
than \cite{OS:schwarzian}; we write $\beta(\rho_1,\rho_2)$ for the
$(2,0)$ part of what Osgood and Stowe call $B_{\rho_1^2}(\log(\rho_2 /
\rho_1))$.}.  For Riemann surfaces, the definition above seems more
natural.

Using this definition, we can compute the Schwarzian tensor in
local coordinates for a pair of conformal metrics:
\begin{equation}
\beta(\rho_1,\rho_2) = \left [ (\sigma_2 - \sigma_1)_{zz}  -
    \left.(\sigma_2)_z\right.^2 + \left . (\sigma_1)_z\right.^2
    \right ] dz^2 \; \text{ where } \; \rho_i = e^{\sigma_i} |dz|
\end{equation}

The Schwarzian tensor generalizes $S(f)$ in the following sense: If
$\Omega \subset \C$ and $\rho$ is the pullback of the Euclidean metric
$|dz|^2$ of $\C$ under a holomorphic map $f : \Omega \to \C$, then
\begin{equation}
\label{eqn:schwarzian-generalization}
\beta(|dz|, \rho) = \beta(|dz|, f^*(|dz|)) = \half S(f).
\end{equation}
Generalizing the cocycle property of the Schwarzian derivative, the
Schwarzian tensors associated to a triple of
conformal metrics $(\rho_1,\rho_2,\rho_3)$ satisfy
\begin{equation}
\label{eqn:schwarzian-cocycle-property}
\beta(\rho_1,\rho_3) = \beta(\rho_1,\rho_2) + \beta(\rho_2,\rho_3).
\end{equation}
Note that $\beta(\rho,\rho)=0$ for any conformal metric $\rho$, so we
also have the antisymmetry relationship:
\begin{equation}
\label{eqn:schwarzian-antisymmetry-property}
\beta(\rho_1,\rho_2) = -\beta(\rho_2,\rho_1)
\end{equation}
Finally, the Schwarzian tensor is functorial with respect to conformal
maps, i.e.
\begin{equation}
\label{eqn:schwarzian-naturality-property}
\beta(f^*(\rho_1),f^*(\rho_2)) = f^*(\beta(\rho_1,\rho_2)),
\end{equation}
where $f$ is a conformal map between domains on Riemann surfaces, and
$\rho_1,\rho_2$ are conformal metrics on the target of $f$.

For a domain $\Omega \subset \C$, we will say a conformal metric is
\emph{M\"obius flat} if its Schwarzian tensor relative to the
Euclidean metric vanishes.  By the cocycle formula, this property is
invariant under pullback by M\"obius transformations.

\begin{lem}[Osgood and Stowe \cite{OS:schwarzian}]
\label{lem:mobius-flat-metrics}
After pulling back by a M\"obius transformation and multiplying by a
positive constant, a M\"obius flat metric on a domain $\Omega \subset
\C$ can be transformed to the restriction of exactly one of the
following examples:
\begin{rmenumerate}
\item The standard Euclidean metric of $\C$,
\label{item:euclidean}
\item The spherical metric of $\CP^1 \simeq S^2 \subset \R^3$,
\label{item:spherical}
\item The hyperbolic metric on a round disk $D \subset \CP^1$,
\label{item:hyperbolic}
\end{rmenumerate}
\end{lem}
We call \ref{item:euclidean}-\ref{item:hyperbolic} the
\emph{Euclidean}, \emph{spherical}, and \emph{hyperbolic} cases,
respectively.  By \eqref{eqn:schwarzian-cocycle-property}, the
property of being M\"obius flat is also equivalent to having vanishing
Schwarzian tensor relative to any other M\"obius flat metric.

While all M\"obius-flat metrics on domains in $\Hat{\C}$ have constant
curvature, the converse is not true.  In fact, a conformal metric has
constant curvature if and only if its Schwarzian relative to a
M\"obius-flat metric is \emph{holomorphic}
\cite[Thm.~1]{OS:schwarzian-of-harmonic-maps}.  For example, we can
calculate the Schwarzian tensor of the (unique up to scale) complete
Euclidean metric $|z^{-1} dz|$ on $\C^*$ relative to the Euclidean
metric of $\C$:
\begin{equation}
\label{eqn:schwarzian-of-cstar-metric}
\beta(|dz|, |z^{-1}\:dz| ) = \frac{1}{4}\frac{dz^2}{z^2}
\end{equation}
Note that everything except the constant $\tfrac{1}{4}$ in
\eqref{eqn:schwarzian-of-cstar-metric} can be derived by symmetry
considerations; the constant itself is determined by calculation.

\section{Decomposition of the Schwarzian}
\label{sec:decomposition-of-the-schwarzian}

In this section we show that the Schwarzian derivative of the
developing map of a grafted surface can be understood geometrically in
terms of the grafting lamination using the Schwarzian tensor of Osgood
and Stowe.  While we restrict attention here to $\lambda \in \ML(S)$
supported on a union of simple closed geodesics, this condition will
be eliminated in the proof of the main theorem by a continuity
argument.  We begin with a brief discussion of the motivation.

When the grafting lamination is a simple closed hyperbolic geodesic,
the restriction of the developing map to the grafted part has a simple
form: it is the composition of a uniformizing map $\Tilde{A} \to \H$
on the (universal cover of) the grafting cylinder $A$ and a map of the
form $z \mapsto z^\alpha$, where $\alpha$ is determined by the measure
on the geodesic.  Since the Schwarzian derivative of a univalent map
is bounded (by a theorem of Nehari), one can use the cocycle property
\eqref{eqn:schwarzian-cocycle} to determine the Schwarzian derivative
of the developing map up to a bounded error.  However, this approach
gives a bound that depends on the homotopy class of the closed
geodesic in an essential way, and offers little hope of an extension
to more general measured laminations.

Rather than expressing part of the developing map as a composition,
the generalized cocycle property
\eqref{eqn:schwarzian-cocycle-property} of the Schwarzian
\emph{tensor} suggests that we look for a conformal metric on the
entire surface $X$ that interpolates between the hyperbolic metric and
the pullback of a spherical metric on $\Hat{\C}$ by the developing
map.  It turns out that the Thurston metric has the right properties
to give a uniform estimate (as we will see in
\S\ref{sec:schwarzian-of-thurston-metric}).

\begin{thm}[Schwarzian decomposition]
\label{thm:schwarzian-decomposition}
Let $X(\lambda) \in P(X)$ be the projective structure on $X$ with
grafting lamination $\lambda \in \ML(S)$, and suppose that $\lambda$
is supported on a union of simple closed geodesics.  Let
$\phi_T(\lambda) \in Q(X)$ be the Schwarzian derivative of its
developing map, and $\Phi(\lambda)$ the grafting differential (which
is not holomorphic).  Then
\begin{equation*}
2 \phi_T(\lambda) = 4 \beta(\rho_0,\rho_\lambda) - \Phi(\lambda).
\end{equation*}
\end{thm}

\begin{proof}
\qedlabel{Theorem \ref{thm:schwarzian-decomposition}}
As the argument is essentially local, we suppress the distinction
between metrics and differentials on $X$ and their lifts to
equivariant objects on $\tilde{X}$.

Let $\rho_{\Hat{\C}}$ be a M\"obius-flat metric on $\Hat{\C}$ (e.g. a
spherical metric).  Using \eqref{eqn:schwarzian-generalization} we have
\begin{equation*}
\phi_T(\lambda) = 2 \beta(\rho_0,f^*\rho_{\Hat{\C}})
\end{equation*}
where $f : \tilde{X} \to \Hat{\C}$ is the developing map.  By the
cocycle property of the Schwarzian tensor,
\begin{equation*}
\beta(\rho_0,f^*\rho_{\Hat{\C}}) = \beta(\rho_0,\rho_\lambda) + \beta(\rho_\lambda,f^*\rho_{\Hat{\C}}).
\end{equation*}
Since the grafting differential is defined as $\Phi(\lambda) = 4
\Phi(\kappa)$, it suffices to show
\begin{equation}
\label{eqn:hopf-of-collapse-is-schwarzian}
\Phi(\kappa) = - \beta(\rho_\lambda,f^*\rho_{\Hat{\C}})
\end{equation}
almost everywhere on $X$, where $\kappa: X \to Y_\lambda = \pr_\lambda
X$ is the collapsing map.

To prove \eqref{eqn:hopf-of-collapse-is-schwarzian}, recall that when
$\lambda$ is supported on a union of simple closed geodesics, the
collapsing map and Thurston metric are based on two local models
\cite[\S 2]{Tanigawa:harmonic-maps}:
\begin{enumerate}
\item In the
hyperbolic part $X_{-1}$, the collapsing map is an isometry, and
\begin{equation*}
\left. \Phi(\kappa)\right|_{X_{-1}}
= 0.
\end{equation*}
On the other hand, the Thurston metric is the pullback by the
developing map of the hyperbolic metric on a round disk in $\Hat{\C}$.
Thus in a neighborhood of a point in the hyperbolic part,
\begin{equation*}
\beta(\rho_\lambda,f^*\rho_{\Hat{\C}})
= \beta(f^*\rho_D,f^*\rho_{\Hat{\C}}) \
= f^* \beta(\rho_D,\rho_{\Hat{\C}}) = 0
\end{equation*}
where $\rho_D$ is the hyperbolic metric on a round disk $D \subset
\Hat{\C}$.  Here we have used the naturality property
\eqref{eqn:schwarzian-naturality-property} of $\beta$ and the fact
that both $\rho_D$ and $\rho_{\Hat{\C}}$ are M\"obius flat.

\item In the Euclidean (grafted) part $X_0$, which is a union of
  cylinders, the collapsing map is locally modeled on the projection
  of $\H$ to $i \R$ by $z \mapsto i |z|$, which has Hopf differential
\begin{equation*}
\Phi(\kappa) = \frac{1}{4} \frac{dz^2}{z^2}.
\end{equation*}
In the same coordinates, the Thurston metric is the pullback
  of the cylindrical metric $|dz|/|z|$ on $\C^*$, so near a point in
  the grafted part,
\begin{equation*}
\begin{split}
\beta(\rho_\lambda,f^*\rho_{\Hat{\C}}) 
&= \beta(f^*\rho_{\C^*},f^*\rho_{\Hat{\C}}) 
= f^* \beta(\rho_{\C^*},\rho_{\Hat{\C}})\\
&= -f^* \beta(\rho_{\Hat{\C}},\rho_{\C^*})
= -\frac{1}{4} \frac{dz^2}{z^2}
\end{split}
\end{equation*}
where in the last line we have used
\eqref{eqn:schwarzian-of-cstar-metric} and the fact that
$\rho_{\Hat{\C}}$ is M\"obius flat.

\end{enumerate}
Thus $\Phi(\kappa:X\to Y_\lambda)$ and
$-\beta(\rho_\lambda,f^*\rho_{\Hat{\C}})$ are equal in $X_{0}$ and
$X_{-1}$, hence a.e. on $X$, which is
\eqref{eqn:hopf-of-collapse-is-schwarzian}, and the theorem follows.
\end{proof}

In light of Theorem \ref{thm:schwarzian-decomposition}, the remaining
obstacle to a geometric understanding of the Schwarzian derivative of
the developing map is the Schwarzian tensor
$\beta(\rho_0,\rho_\lambda)$ of the Thurston metric relative to the
hyperbolic metric.  After studying the Thurston metric in more detail,
we will determine a bound for its Schwarzian tensor in
\S\ref{sec:schwarzian-of-thurston-metric}.

\section{NPC conformal metrics}
\label{sec:npc-metrics}

Some of the properties of the Thurston metric on a grafted surface
that we will use in the proof of the main theorem can be attributed to
the fact that it is nonpositively curved (NPC).  We devote this
section and the next to a separate discussion of such metrics and the
additional regularity properties they enjoy compared to general
conformal metrics on a Riemann surface.

Consider a geodesic metric space $(M,d)$, i.e. a metric space in which
the distance $d(x,y)$ is then length of some path joining $x$ and $y$.
We say $(M,d)$ is nonpositively curved (NPC) space if all of its
geodesic triangles are ``thinner'' than triangles in the plane with
the same edge lengths (see \cite{ABN} for details on this and
equivalent definitions).  A space with this property is also called
$\cat(0)$.  The definition of an NPC space actually implies that it is
simply connected, but we will also say that a metric manifold $(M,d)$
is NPC if the triangle condition is satisfied in its universal cover
$(\tilde{M},\tilde{d})$.

An NPC metric on a surface $S$ naturally induces a conformal structure:

\begin{thm}
\label{thm:conformal-representation-of-npc-surface}
Let $S$ be a compact surface and $d(\cdot,\cdot)$ an NPC metric
compatible with the topology of $S$.  Then there is a unique Riemann
surface $X \in \T(S)$ and a conformal metric $\rho$ on $X$ inducing
$d$, i.e. such that
\begin{equation*}
 d(x,y) = \inf \: \left \{ \left. \int_\gamma \rho(z) \: |dz| \; \right | \; \gamma : ([0,1],0,1) \to
  (X,x,y) \right \}.
\end{equation*}
\end{thm}

Theorem \ref{thm:conformal-representation-of-npc-surface} is
essentially due to Re\v{s}etnyak
(see \cite{Reshetnyak:isothermal-coordinate},
\cite{Reshetnyak:nonexpanding-maps}, and the recent survey
\cite{Reshetnyak:conformal-representation}) who shows that for an NPC
metric on a two-dimensional manifold there is a local conformal
homeomorphism to the disk $\DD \subset \C$.  Here ``conformal'' must be
interpreted in terms of the preservation of angles between curves,
which are defined in an NPC space using the distance function.  A
theorem of Huber implies that under
these circumstances there is a conformal metric on the disk that gives
$\rho$ as above \cite{Huber:alexandrov-spaces}.

Another proof of Theorem
\ref{thm:conformal-representation-of-npc-surface} is given in
\cite{Mese:singular-spaces} using the Korevaar-Schoen theory of
harmonic maps to metric spaces; here the metric $\rho$ is obtained as
the pullback metric tensor for a conformal harmonic map from a domain
in $\C$ to a domain in the NPC surface $(S,d)$.  This analysis leads
naturally to more detailed regularity and nondegeneracy properties of
$\rho$:

\begin{thm}[Mese \cite{Mese:singular-spaces}]
There is a one-to-one correspondence between NPC metrics
$d(\cdot,\cdot)$ on $S$ and the distance functions arising from pairs
$(X,\rho)$ with $X \in \T(S)$ and $\rho$ a conformal metric on $X$
such that $\rho \in W^{1,2}_{\loc}(X)$ and $\log \rho(z)$ is weakly
subharmonic.  Furthermore, if $\rho$ is such a conformal metric, then
$\rho(z) > 0$ almost everywhere.
\end{thm}

The subharmonicity of $\log \rho(z)$ reflects the condition of
nonpositive curvature; indeed if $\rho$ is a smooth, nondegenerate
conformal metric then its Gaussian curvature at a point $z$ is
\begin{equation}
\label{eqn:gaussian-curvature}
K_\rho(z) = -\frac{\Delta \log \rho(z)}{\rho(z)^2}
\end{equation}
and therefore
\begin{equation*}
K_\rho(z) \leq 0 \; \Leftrightarrow \; \Delta \log \rho(z) \geq 0.
\end{equation*}
For more general NPC conformal metrics on $X$, there is no direct
analogue of the Gaussian curvature function, but there is a curvature
measure $\Omega_\rho$ with local expression
\begin{equation}
\label{eqn:curvature-laplacian}
\Omega_{\rho} = - \Delta \log \rho
\end{equation}
By approximation one can show that the
Gauss-Bonnet theorem holds in this context, i.e. the total measure of
$\Omega_\rho$ is $2 \pi \chi(S)$.

So finally we define $\npc(X)$ to be the set of conformal metrics
$\rho \in L^2(S_{\half,\half}(X))$ such that the induced distance
function $d_\rho(\cdot,\cdot)$ makes $X$ into an NPC metric space.
From the preceding discussion, $\rho(z)$ is then $W^{1,2}_{\loc}$,
almost everywhere positive, $\log \rho(z)$ is subharmonic in a
conformal coordinate chart, and the area $dA_\rho$ and curvature
measures $\Omega_\rho$ are finite.  Conversely, any NPC metric on $S$
whose associated conformal structure (as in Theorem
\ref{thm:conformal-representation-of-npc-surface}) is $X$ gives rise
to such a conformal metric.

\section{Regularity and compactness for NPC metrics}
\label{sec:regularity-and-compactness}

We now compare the various NPC metrics on a fixed compact Riemann
surface $X$.  Let $\npc_a(X) \subset \npc(X)$ denote the set of NPC
metrics on $X$ with total area $2 \pi |\chi| a$ (i.e.~$a$ times the
hyperbolic area of $X$).  For any $t>0$, we have $\rho \in \npc_a(X)$
if and only if $t\rho \in \npc_{t^2a}(X)$, so we may as well consider
only metrics of a fixed area, e.g. $\rho \in \npc_1(X)$.

\begin{thm}[{Mese \cite[Thm.~29]{Mese:compactness}}]
\label{thm:czero-compactness}
The set of distance functions $\{ d_\rho \: | \: \rho \in \npc_1(X)
\}$ is compact in the topology of uniform convergence.  In particular,
it is closed: a uniform limit of such distance functions is the
distance function of an NPC metric on $X$.
\end{thm}

\begin{rem}
In \cite{Mese:compactness}, Mese establishes this compactness result
in the greater generality of metrics with a positive upper bound on
curvature (i.e. $\cat(k)$, rather than $\cat(0)$) using the theory of
harmonic maps to such metric spaces.
\end{rem}

Conformal metrics $\rho_1$ and $\rho_2$ with uniformly close distance
functions can differ significantly on a small scale; for example, such
an estimate does not give any control on the modulus of continuity of
the map of metric spaces $\id:(X,\rho_1) \to (X,\rho_2)$.  Since this
is precisely the kind of control we will need for application to
grafting, we now investigate the local properties of NPC metrics.
Since it is smooth and uniquely determined, the hyperbolic metric
$\rho_0 \in \npc_1(X)$ is a good basis for comparison of regularity of
NPC metrics on $X$.

\begin{thm}
\label{thm:lipschitz-holder}
For $X \in \T(S)$ with hyperbolic metric $\rho_0 \in \npc_1(X)$ we
have:
\begin{rmenumerate}
\item For all $\rho \in \npc_1(X)$, the map $\id : (X,\rho_0) \to
  (X,\rho)$ is Lipschitz with constant depending only on $X$, i.e.
 \label{item:lipschitz}
\begin{equation*} \| \rho / \rho_0 \|_\infty  \leq C(X). \end{equation*}

\item For all $\rho_1,\rho_2 \in \npc_1(X)$, the map $\id : (X,\rho_1)
  \to (X,\rho_2)$ is bi-H\"older, with exponent and an upper bound on
  the H\"older norm depending only on $X$.
\label{item:holder}
\end{rmenumerate}
\end{thm}

\begin{rem}
Modulo the values of exponents and Lipschitz and H\"older norms,
Theorem \ref{thm:lipschitz-holder} is the best kind of estimate that
could hold for a general NPC metric of fixed area: When $\rho =
|\phi|^{\half}$ for a holomorphic quadratic differential $\phi$, the
map $(X,\rho) \to (X,\rho_0)$ is H\"older but not Lipschitz near the
zeros of $\phi$.
\end{rem}

\begin{proof}[Proof (of Theorem \ref{thm:lipschitz-holder})]\mbox{}

\noindent \ref{item:lipschitz} Consider the map $\id : (X,\rho_0) \to
  (X,\rho)$, which is conformal and thus harmonic.  Its energy is the
  area of the image, $2 \pi |\chi|$, thus by Korevaar and Schoen's
  regularity theorem for harmonic maps to NPC spaces, the Lipschitz
  constant of the map is bounded above by a constant that depends only
  on the hyperbolic metric $\rho_0$, i.e. on the conformal structure
  $X$ \cite{KS:sobolev-spaces}.

\medskip
\noindent \ref{item:holder} In light of \ref{item:lipschitz}, it
suffices to give a lower bound on the $(X,\rho)$-distance in terms of
the $(X,\rho_0)$-distance for all $\rho \in \npc_1(X)$.  That is, for
$x,y \in X$ sufficiently close, we must show that
\begin{equation}
\label{eqn:distance-goal}
d_\rho(x,y) \geq C d_{\rho_0}(x,y)^K
\end{equation}
for some $C,K > 0$ that depend on $X$.

To establish \eqref{eqn:distance-goal}, we will use the fact that
$\log(\rho)$ is subharmonic (by the NPC condition), which limits the
size of the set where $\log(\rho)$ is large and negative (making
$\rho$ close to zero).  The distance estimate comes from an effective
form of the classical fact that the set where a subharmonic function
takes the value $-\infty$ has zero Hausdorff dimension; this will
prevent $\rho$ from being too small on a significant fraction of a
geodesic segment.

First we make the problem local by covering $X$ with disks where the
metric $\rho$ is everywhere bounded above and somewhere bounded below.
Specifically, we cover $X$ by $N$ hyperbolic disks of a fixed radius
$R$, each parameterized by the unit disk $\DD \subset \C$.  Here $N$
and $R$ depend only on $X$.  We can also find $r<1$ depending on $X$
such that the images of $\DD_{r} = \{ |z| < r \}$ under the $N$ chart
maps still cover $X$.

For such a covering by disks, there are constants $M > m > -\infty$,
also depending only on $X$, such that on each such disk the conformal
density $\rho(z)$ satisfies
\begin{equation}
\label{eqn:rho-bounded-above-below}
\begin{split}
\log \rho(z) < M \text{ for all } z \in \DD\\
\log \rho(z_0) > m \text{ for some } z_0 \in \DD_{r}
\end{split}
\end{equation}
The existence of the upper bound $M$ follows from the Lipschitz
estimate in part \ref{item:lipschitz} above.  As for $m$, if such a
bound did not exist, one could find a sequence $\rho_i \in \npc_1(X)$
converging uniformly to $0$ on $\DD_{r}$.  In particular the
associated distance functions $d_{\rho_i}$ could not accumulate on the
distance function of any NPC metric $d_\infty$, contradicting Theorem
\ref{thm:czero-compactness}.

We now show that for a conformal density function $\rho$ of an NPC
metric on the unit disk satisfying \eqref{eqn:rho-bounded-above-below}
and any $x,y \in \DD_r$ sufficiently close,
\begin{equation}
\label{eqn:distance-power}
d_\rho(x,y) \geq C |x-y|^K
\end{equation}
where $C$ and $K$ depend only on $X$.  Since the hyperbolic metric
$\rho_0$ is smooth and comparable to the Euclidean metric (by constants
depending on $X$ and the covering by disks), part \ref{item:holder}
of the theorem then follows with a different constant $C$.

Fix $r'$ such that $r < r' < 1$.  For $x,y \in \DD_r$ close enough,
there is a minimizing $\rho$-geodesic segment $\gamma \subset
\DD_{r'}$ joining them, because by Theorem \ref{thm:czero-compactness}
the $\rho$-distance from $\partial \DD_{r}$ to $\partial \DD_{r'}$ is
bounded below for all $\rho \in \npc_1(X)$.

We now use a result of Brudnyi on the level sets of subharmonic
functions:

\begin{thm}[{Brudnyi \cite[Prop.\ 1]{Brudnyi:subharmonic-bmo}}]
\label{thm:brudnyi-subharmonic-bound}
Let $u$ be a subharmonic function on $\DD$ such that
\begin{equation*} \sup_{\DD} u < M \; \text{ and } \; \sup_{\DD_r} u > m \end{equation*}
for some $r<1$ and $-\infty < m < M < \infty$.  Then for each
$\epsilon>0$ and $d>0$ there is a countable family of disks
$\DD(z_i,r_i)$ with
\begin{equation*} \sum_i r_i^d < \frac{(2\epsilon)^d}{d} \end{equation*}
and such that
\begin{equation*} u(z) \geq m -  C (M - m) \left ( 1 + \log\left(\frac{1}{\epsilon}\right) \right ) \end{equation*}
for all $z \in \DD_r$ outside these disks.  Here $C=C(r) > 0$ depends
only on $r$.
\end{thm}

Applying Theorem
\ref{thm:brudnyi-subharmonic-bound}, there is a collection of disks in
$\DD$ with radii $r_i$ satisfying
\begin{equation*} \sum_i r_i \leq \frac{|x-y|}{2} \end{equation*}
and such that on the complement of these disks in $\DD_{r'}$, the conformal factor
$\log \rho(z)$ satisfies
\begin{equation}
\label{eqn:log-rho-bound}
\log \rho(z) > m - C (M - m) \left ( 1 + \log \frac{4}{|x-y|} \right )
\end{equation}
Since the path $\gamma$ connects $x$ and $y$, the
exclusion of these disks leaves a subset $\gamma_0 \subset \gamma$ of Euclidean length
at least $|x-y| / 2$ where the bound \eqref{eqn:log-rho-bound} is
satisfied.  In particular,
\begin{equation}
\begin{split}
d_\rho(x,y) &= \int_\gamma \rho(z) \: |dz| \geq \int_{\gamma_0} \rho(z) \: |dz|\\
& \geq \exp\left( \inf_{\gamma_0} \log \rho(z) \right ) \int_{\gamma_0} \: |dz|\\
& \geq \exp\left(C_1 - C_2\log(4 / |x-y|)\right ) \: \frac{|x-y|}{2} \\
& \geq C_3 \: |x-y|^K
\end{split}
\end{equation}
where $C_i$ are unspecified constants that depend only on $X$.  This
is the desired bound from \eqref{eqn:distance-power}, so we have
proved \ref{item:holder}.

\qedlabel{Theorem \ref{thm:lipschitz-holder}}
\end{proof}

\section{Curvature of the Thurston metric}
\label{sec:thurston-metric}

\begin{figure}
\begin{center}
\includegraphics[width=13cm]{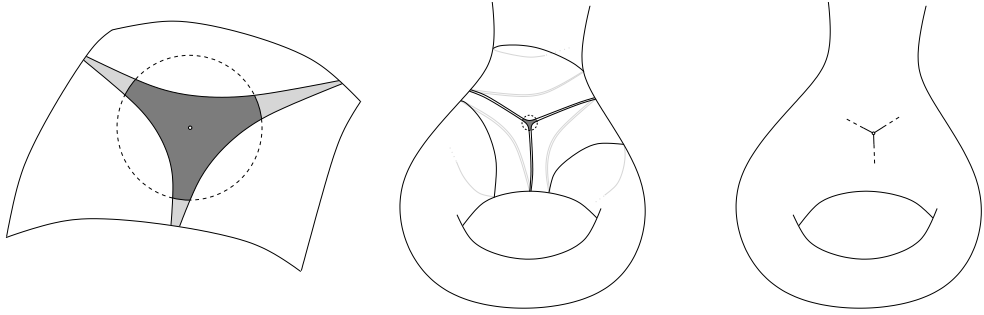}
\caption{Scaling the Thurston metric (left) to have unit area makes
  its curvature concentrate near a finite set of points (center) in a
  manner reminiscent of a quadratic differential metric (right).}
\label{fig:thurston-metric}
\end{center}
\end{figure}

We now begin to study the geometry of the Thurston metric.  The idea
is that for large grafting laminations $\lambda$, the Thurston metric
looks a lot like the singular flat metric coming from a holomorphic
quadratic differential.  This is because the $\rho_\lambda$-area of
the hyperbolic part $X_{-1}(\lambda)$ is fixed (and equal to the
hyperbolic area of $Y_\lambda$), while that of the Euclidean part
$X_0(\lambda)$ grows with $\lambda$.  In fact, we have (cf. \cite[\S
3]{Tanigawa:harmonic-maps})
\begin{equation*}
\area(X_0(\lambda), \rho_\lambda) = \ell(\lambda,Y_\lambda),
\end{equation*}
while by Theorem \ref{thm:grafting-energy-close}, 
\begin{equation*}
\ell(\lambda,Y_\lambda) \approx E(\lambda,X) \to \infty \text{ as } \lambda \to
\infty,
\end{equation*}
where $E(\lambda,X)$ is the extremal length of $\lambda$ on $X$.  Thus
if we rescale the Thurston metric to have constant area, all of its
curvature is concentrated in a very small part of the surface, which
is reminiscent of the conical singularities of a quadratic
differential metric (see Figure \ref{fig:thurston-metric}).

Since the area of $(X,\rho_\lambda)$ is approximately $E(\lambda,X)$
(i.e. up to a bounded additive constant), the ratio $\rho_\lambda /
\rho_0$ is approximately $E(\lambda,X)^{\half}$ in an average
sense.  Since $\rho_\lambda /\!\area(X,\rho_\lambda)^{\half} \in
\npc_1(X)$, part \ref{item:lipschitz} of Theorem
\ref{thm:lipschitz-holder} provides an upper bound of the same order, i.e.
\begin{equation}
\label{eqn:pointwise-bound-for-thurston}
\frac{\rho_\lambda(x)}{\rho_0(x)} \leq C(X)\left (1 +
E(\lambda,X)^{\half} \right )
\end{equation}

The next theorem quantifies the sense in which the curvature of
$\rho_\lambda$ becomes concentrated for large $\lambda$: A finite set
of small hyperbolic disks on $X$ suffices to cover all but an
exponentially small part of $X_{-1}(\lambda)$, measured with respect
to the background metric $\rho_0$.

\begin{thm}[curvature concentration]
\label{thm:curvature-concentration}
For any $\epsilon > 0$, $X \in \T(S)$ and $\lambda \in \ML(S)$, let
$X_{-1}(\lambda)$ denote the subset of $X$ where the Thurston metric
$\rho_\lambda$ is hyperbolic.  Then there are $N = 2 |\chi|$
points $x_1, \ldots, x_N \in X$ such that
\begin{equation*} \area\left( X_{-1}(\lambda) - B, \rho_0 \right) \leq  C \exp(-\alpha E(\lambda,X)^{\half})),\end{equation*}
 where $B = \bigcup_i B_\epsilon(x_i)$, $B_\epsilon(x)$ is the
hyperbolic ball of radius $\epsilon$ centered at $x$, and the
constants $C$ and $\alpha$ depend on $\epsilon$ and $X$ (but not on
$\lambda$).
\end{thm}

The curvature concentration phenomenon described by Theorem
\ref{thm:curvature-concentration} arises from a simple geometric
property of ideal triangles.

\begin{lem}[ideal triangles]
\label{lem:ideal-triangles}
Let $T \subset \H^2$ be a hyperbolic ideal triangle (a region
bounded by three pairwise asymptotic geodesics), and $x \in T$.  Then
\begin{equation*} \area(T - B_r(x)) \leq C e^{-r} \end{equation*}
where $C$ depends only on $d(x,\partial T)$. 
\end{lem}

\begin{proof}
This is an exercise in hyperbolic geometry; an explicit calculation
shows that this is true for the center of $T$, and then the full
statement follows since for each $\epsilon > 0$, $(T -
N_\epsilon(\partial T))$ is compact.  \qedlabel{Lemma
\ref{lem:ideal-triangles}}
\end{proof}

\begin{proof}[Proof of Theorem \ref{thm:curvature-concentration}.]
Fix $\epsilon > 0$.  Every geodesic lamination on a hyperbolic surface
can be enlarged (non-uniquely) to a geodesic lamination whose
complement is a union of $N = 2 |\chi|$ hyperbolic ideal triangles.
Given $\lambda \in \ML(S)$, let $\tau$ be such an enlargement of the
supporting geodesic lamination of $\lambda$ on the hyperbolic surface
$Y_\lambda = \pr_\lambda X$.

Let $T_1, \ldots, T_N$ be the ideal triangles comprising $Y_\lambda -
\tau$, and $t_1, \ldots, t_N \in Y_\lambda$ points in the thick parts
of these triangles (e.g.~their centers).  Since $X = \gr_\lambda
Y_\lambda$, each of the triangles $T_i$ naturally includes into $X$,
and this inclusion is an isometry for the Thurston metric.  Let $x_i
\in X$ denote the point corresponding to $t_i \in T_i$ by this
inclusion.  The images of $T_i$ in $X$ cover all of $X_{-1}(\lambda)$
except a null set consisting of the finitely many
$\rho_\lambda$-geodesics added to $\supp(\lambda)$ to obtain $\tau$.

Since $\rho_\lambda /\!\area(X,\rho_\lambda)^{\half} \in \npc_1(X)$,
part \ref{item:holder} of Theorem \ref{thm:lipschitz-holder} implies
that the metrics $\rho_0$ and $\rho_\lambda
/\!\area(X,\rho_\lambda)^{\half}$ are H\"older equivalent.  Thus for
some $k > 0$, the hyperbolic disk $B_\epsilon(x_i)$ contains a
$\rho_\lambda$-disk of radius at least $C(X) E(\lambda,X)^{\half}
\epsilon^k$, where $C$ and $k$ depend only on $X$.  Here we have used
the fact that $\area(X,\rho_\lambda) \approx E(\lambda,X)$ and
$\area(X,\rho_\lambda)$ is bounded below (by $\area(X,\rho_0) = 2 \pi
|\chi|$).

Applying Lemma \ref{lem:ideal-triangles} to $T_i$, we see that
$B_\epsilon(x_i)$ covers all of the image of $T_i$ in $X$ except for a
region of $\rho_\lambda$-area at most $C_1(X) \exp(-C_2(X)
E(\lambda,X)^{\half} \epsilon^k)$.
Since $\rho_\lambda \geq \rho_0$, the same upper bound holds for the
hyperbolic area, and applying this to each of the $N$ triangles we obtain:
\begin{equation*}
\begin{split}
\area(X_{-1}(\lambda) - B,\rho_0) &\leq \area(X_{-1}(\lambda) -
B,\rho_\lambda)\\
&= \area\left( \bigcup_iT_i - \bigcup_iB_\epsilon(x_i),
\rho_\lambda \right)\\
&\leq N C_1(X) \exp(-C_2(X) E(\lambda,X)^\half \epsilon^k )
\end{split}
\end{equation*}
Taking $\alpha = C_2(X) \epsilon^k$, Theorem
\ref{thm:curvature-concentration} follows.

\qedlabel{Theorem \ref{thm:curvature-concentration}}
\end{proof}

\section{The Schwarzian tensor of the Thurston metric}
\label{sec:schwarzian-of-thurston-metric}

In this section we find an upper bound for the norm of the Schwarzian
tensor of the Thurston metric restricted to a subset of $X$.  As in
Theorem \ref{thm:curvature-concentration}, we could take this subset to
be the complement of finitely many small hyperbolic balls, but it is
technically simpler to work with a fixed domain, so we use:

\begin{lem}
\label{lem:delta}
For any $X \in \T(S)$ there exist $\epsilon_0(X), \delta(X) > 0$ such
that if $B \subset X$ is the union of $N = 2|\chi|$ hyperbolic balls
of radius $\epsilon \leq \epsilon_0(X)$, then $(X-B)$ contains an
embedded hyperbolic ball of radius $\delta(X)$.
\end{lem}

\begin{proof}
Choose $\epsilon_0$ small enough so that $2 \epsilon_0$ is less than
the hyperbolic injectivity radius of $X$ and so that $2N$ balls of
radius $2 \epsilon_0$ cannot cover $X$ (say, by area considerations).
Then for any union of $N$ balls of radius $\epsilon < \epsilon_0$,
there is a point in $X$ whose distance from all of them is at
least $\epsilon_0$.  Thus $(X-B)$ contains a hyperbolic ball $D$ of
radius $\epsilon_0$, which is necessarily embedded, for any $B$ as in
the statement of the Lemma.  We set $\delta = \epsilon_0$.
\qedlabel{Lemma \ref{lem:delta}}
\end{proof}

Now we continue our study of the Thurston metric with a series of
analytic results that apply Theorem \ref{thm:curvature-concentration}
(curvature concentration).  Recall from \eqref{eqn:sigma} that
$\sigma(\rho_1,\rho_2) = \log(\rho_2/\rho_1)$.
\begin{lem}
\label{lem:lp-norm-laplacian}
Fix $X \in \T(S)$ and $p < \infty$.  For each $\lambda \in \ML(S)$
there is a hyperbolic ball $D \subset X$ of radius $\delta(X)$ such
that
\begin{equation*}
\| \Delta_{\rho_0} \sigma(\rho_0, \rho_\lambda) \|_{L^p(D,\rho_0)} < C(p,X).
\end{equation*}
\end{lem}

\begin{proof}
Using the formula for the Gaussian curvature of a conformal metric
\eqref{eqn:gaussian-curvature}, we have
\begin{equation*}
\Delta_{\rho_0} \sigma(\rho_\lambda,\rho_0) = -\frac{1}{\rho_0^2}
\Delta \log \rho_0 + \frac{1}{\rho_0^2} \Delta \log \rho_\lambda = K_{\rho_0} - \frac{\rho_\lambda^2}{\rho_0^2} K_{\rho_\lambda}.
\end{equation*}
Note that the Gaussian curvature of the Thurston metric exists almost
everywhere because its conformal factor is $C^{1,1}$
\cite{KP:canonical-metric}.  Since $K_{\rho_0} \equiv -1$ and the
$\rho_0$-area of $X$ is $2 \pi |\chi|$ we have
\begin{equation*}
\| K_{\rho_0}\|_{L^p(X,\rho_0)} = \left ( 2 \pi |\chi| \right )^\frac{1}{p} =
C(p)
\end{equation*}

So we need only establish an $L^p$ bound for the second term.  From
\eqref{eqn:pointwise-bound-for-thurston}, we have
$\rho_\lambda^2/\rho_0^2 \leq C \: ( 1 +  E(\lambda,X))$.  On the other hand,
$|K_{\rho_\lambda}| \leq 1$, and this function is supported in
$X_{-1}(\lambda)$.  By Theorem \ref{thm:curvature-concentration} we
have
\begin{equation*}
\area(X_{-1}(\lambda) - B,\rho_0) \leq C \exp \left (-\alpha
E(\lambda,X)^{\half} \right )
\end{equation*}
Where $C$ and $\alpha$ depend on $X$ and $\epsilon$.  Combining these
estimates, we find
\begin{equation*}
\begin{split}
\| \frac{\rho_\lambda^2}{\rho_0^2} K_{\rho_\lambda} \|_{L^p(X-B,\rho_0)} &\leq 
\left ( \sup \frac{\rho_\lambda^2}{\rho_0^2}\left |
K_{\rho_\lambda}\right | \right) 
\left ( \area\left ( \supp K_{\rho_\lambda} \cap (X - B), \rho_0\right)
\right )^{\frac{1}{p}} \\
&\leq C \left (1 + E(\lambda,X) \right ) \exp \left( - \alpha p^{-1} E(\lambda,X)^{\half}
\right).
\end{split}
\end{equation*}
In particular, $\| \frac{\rho_\lambda^2}{\rho_0^2} K_{\rho_\lambda}
\|_{L^p(X-B,\rho_0)} \to 0$ as $E(\lambda,X) \to \infty$, and we have
a uniform upper bound on the norm depending only on $X$, $\epsilon$,
and $p$.

Taking $\epsilon=\epsilon_0(X)$ and applying Lemma \ref{lem:delta},
the set $(X-B)$ contains a hyperbolic ball $D$ of radius $\delta(X)$,
and restriction to this set only decreases the norm.  So on $D$ we
obtain an upper bound depending on $p$ and $X$.
\qedlabel{Lemma \ref{lem:lp-norm-laplacian}}
\end{proof}

In what follows it will be convenient to normalize the area of the
Thurston metric; let $\Hat{\rho}_\lambda$ denote the positive multiple
of the Thurston metric with the same area as the hyperbolic metric
$\rho_0$ on $X$, i.e.
\begin{equation*}
\Hat{\rho}_\lambda = \left ( \frac{2 \pi |\chi|}{\area(X,\rho_\lambda)}
  \right )^\half \rho_\lambda.
\end{equation*}

\begin{lem}
\label{lem:lp-norm-sigma}
Fix $X \in \T(S)$ and $p < \infty$.  For each $\lambda \in \ML(S)$
there is a hyperbolic ball $D_{\delta/2} \subset X$ of radius
$\delta(X)/2$ such that
\begin{equation*}
\| \sigma(\rho_0, \Hat{\rho}_\lambda) \|_{L^p(D_{\delta/2})} < C(p,X).
\end{equation*}
\end{lem}

\begin{rem}
Using $\Hat{\rho}_\lambda$ instead of $\rho_\lambda$ only changes
$\sigma(\rho_0,\rho) = \log(\rho / \rho_0)$ by a constant, so
Lemma \ref{lem:lp-norm-sigma} says that
$\sigma(\rho_0,\rho_\lambda)$ is close to a constant in an $L^p$
sense.
\end{rem}

\begin{proof}
Since $\Hat{\rho}_\lambda \in \npc_1(X)$, part \ref{item:lipschitz} of
Theorem \ref{thm:lipschitz-holder} implies that
$\sigma(\rho_0,\Hat{\rho}_\lambda) < C_0(X)$.  We will use the bound
on the curvature of $\Hat{\rho}_\lambda$ to turn this global upper
bound into an estimate of the $L^p$ norm in a small disk.

Let $F = C_0 - \sigma(\rho_0,\Hat{\rho}_\lambda)$, so $F$ is a
nonnegative function, and let $D = D_\delta$ be the hyperbolic ball of
radius $\delta$ provided by Lemma \ref{lem:lp-norm-laplacian}.
Then we have
\begin{equation*}
\| \Delta_{\rho_0} F\|_{L^p(D_\delta)} = \| \Delta_{\rho_0} \sigma(\rho_0,\Hat{\rho}_\lambda) \|_{L^p(D_\delta)} \leq C_1(p,X).
\end{equation*}

Let $D_{r}$ denote the hyperbolic ball of radius $r$ concentric with
$D$.  We now want to show that $F$ cannot be large throughout
$D_{\delta/4}$, or equivalently, show that there is a point $z \in
D_{\delta/4}$ where $\Hat{\rho}_\lambda = \exp(C_0 - F) \rho_0$ is not
too close to zero.  But as we saw in the proof of Theorem
\ref{thm:lipschitz-holder}, this follows from the compactness of
distance functions for metrics in $\npc_1(X)$ (Theorem
\ref{thm:czero-compactness}): if $\Hat{\rho}_\lambda$ could be
arbitrarily small throughout a disk of definite hyperbolic radius (in
this case, $\delta(X)/4$), a limiting argument would give a metric in
$\npc_1(X)$ that vanishes on an open set, a contradiction.  Thus
$\inf_{D_{\delta/4}} F \leq C_2(X)$.

Combining these bounds on $F$ and $\Delta F$, we apply the weak
Harnack inequality (see \cite[Thm.~8.18]{GT:second-order-pde}) to obtain
\begin{equation}
\label{eqn:harnack}
\begin{split}
\| F \|_{L^p(D_{\delta/2})} &\leq C_3(\delta(X)) \: \left ( \|
\Delta_{\rho_0} F \|_{L^p(D_{\delta})} + \inf_{D_{\delta/4}} F \right
)\\ &\leq C_3(\delta(X)) \left ( C_1(p,X) + C_2(X) \right )\\
&= C(p,X)
\end{split}
\end{equation}
Since $F = C_0 - \sigma(\rho_0,\Hat{\rho}_\lambda)$, the Lemma
follows by algebra.
\qedlabel{Lemma \ref{lem:lp-norm-sigma}}
\end{proof}

Combining the preceding lemmas and standard elliptic theory, we can
now bound the Schwarzian tensor:

\begin{thm}
\label{thm:schwarzian-of-thurston-metric-bounded}
Fix $X \in \T(S)$.  For each $\lambda \in \ML(S)$ there is a
hyperbolic ball $D_{\delta/4} \subset X$ of radius $\delta(X)/4$ such
that
\begin{equation*}
\| \beta(\rho_0,\rho_\lambda) \|_{L^1(D_{\delta/4})} \leq C(X).
\end{equation*}
\end{thm}

\begin{proof}
Let $D_{\delta/2}$ be as in Lemma \ref{lem:lp-norm-sigma} and let
$D_{\delta/4} \subset D_{\delta/2}$ be the concentric ball of radius
$\delta(X)/4$.  By standard elliptic theory
(e.g. \cite[Thm.~9.11]{GT:second-order-pde}) we have the Sobolev norm
estimate
\begin{equation*}
\| u \|_{W^{2,2}(D_{\delta/4},\rho_0)} \leq C(\delta) \left ( \|
  \Delta_{\rho_0} u \|_{L^2(D_{\delta/2},\rho_0)} + \| u
  \|_{L^2(D_{\delta/2},\rho_0)} \right ),
\end{equation*}
that is, the second derivatives of $u$ are bounded in terms of $u$ and
its Laplacian.

Applying this to $u = \sigma(\rho_0,\Hat{\rho}_\lambda)$, the terms on
the right hand side are bounded by Lemmas
\ref{lem:lp-norm-laplacian} and \ref{lem:lp-norm-sigma}, respectively, giving
\begin{equation}
\label{eqn:w22-bound}
\| \sigma(\rho_0,\Hat{\rho}_\lambda) \|_{W^{2,2}(D_{\delta/4},\rho_0)} \leq C_1(X).
\end{equation}

Since the $W^{2,2}$ norm bounds both the derivative of
$\sigma(\rho_\lambda,\rho_0)$ in $L^2$ and its Hessian in $L^2$ (hence
also $L^1$, by the Cauchy-Schwartz inequality), the definition of the
Schwarzian tensor \eqref{eqn:hessian-definition} gives
\begin{equation*}
\| \beta(\rho_\lambda,\rho_0) \|_{L^1(D_{\delta/4},\rho_0)} < C_2(X)
\| \sigma(\rho_0,\Hat{\rho}_\lambda) \|_{W^{2,2}(D_{\delta/4},\rho_0)},
\end{equation*}
which together with \eqref{eqn:w22-bound} gives the desired bound.
\qedlabel{Theorem \ref{thm:schwarzian-of-thurston-metric-bounded}}
\end{proof}

\section{Proof of the main theorem}
\label{sec:proof}

Now that we have an estimate on the Schwarzian tensor of the Thurston
metric (Theorem \ref{thm:schwarzian-of-thurston-metric-bounded}) and
the decomposition of the Schwarzian derivative of a $\CP^1$ structure
(Theorem \ref{thm:schwarzian-decomposition}), the proof of the main
theorem is straightforward.  We will need a lemma about holomorphic
quadratic differentials in order to extend a bound on a small
hyperbolic ball to one on $X$.

\begin{lem}
\label{lem:extension}
For any $X \in \T(S)$ and $\delta > 0$ there is a constant
$C(\delta,X)$ such that
\begin{equation*}
\| \psi \|_{L^1(X)} \leq C(\delta, X) \| \psi \|_{L^1(D_{\delta})}
\end{equation*}
for all $\psi \in Q(X)$ and any hyperbolic ball $D_{\delta} \subset X$
of radius $\delta$.
\end{lem}

\begin{proof}
The inequality is homogeneous, so we need only prove it for $\psi \in
Q(X)$ with $\| \psi \|_{L^1(X)} = 1$.  The set of all unit-norm
quadratic differentials is compact and equicontinuous.

Suppose on the contrary that there is no such constant $C(\delta,X)$.
Then there is a sequence $\psi_n \in Q(X)$ and $\delta$-balls $D_n
\subset X$ such that $\| \psi_n \|_{L^1(X)} = 1$ and $\| \psi_n
\|_{L^1(D_n)} \to 0$ as $n \to \infty$.  Taking a subsequence we can
assume $\psi_n \to \psi_\infty$ uniformly and $D_n \to D_\infty$,
where we say a sequence of $\delta$-balls converges if their centers
converge.  By uniform convergence of $\psi_n$ we obtain
\begin{equation*}
\|\psi_\infty\|_{L^1(D_\infty)} = \lim_{n \to \infty} \| \psi_n
\|_{L^1(D_n)} = 0.
\end{equation*}
But $\psi_\infty$ is a nonzero holomorphic quadratic differential, which
 vanishes at only finitely many points, so this is a contradiction.
\qedlabel{Lemma \ref{lem:extension}}
\end{proof}

\begin{proof}[Proof (of Theorem \ref{thm:schwarzian-nearly-strebel}).]
First suppose that $\lambda$ is supported on a union of simple closed
geodesics.

Comparing both $\phi_F(\lambda)$ and $\phi_T(\lambda)$ to the grafting
differential $\Phi(\lambda)$, and applying Theorem
\ref{thm:schwarzian-decomposition}, we have
\begin{equation*}
\begin{split}
2 \phi_T(\lambda) + \phi_F(\lambda) &= \Bigl ( 2 \phi_T(\lambda) +
  \Phi(\lambda) \Bigr ) \: - \: \Bigl ( \Phi(\lambda) -
  \phi_F(\lambda) \Bigr )\\
&= 4 \beta(\rho_\lambda,\rho_0) - \Bigl (
  \Phi(\lambda) - \phi_F(\lambda) \Bigr )
\end{split}
\end{equation*}
Taking the $L^1$ norm, we apply Theorem
\ref{thm:schwarzian-of-thurston-metric-bounded} to the first term and
Theorem \ref{thm:grafting-differentials-close} to the second, giving
\begin{equation}
\label{eqn:partialbound}
\| 2 \phi_T(\lambda) + \phi_F(\lambda) \|_{L^1(D_{\delta/4})} \leq 4
  C_1(X) + C_2 \left ( 1 + E(\lambda,X)^\half \right)
\end{equation}
where $D_{\delta/4}$ is a hyperbolic ball of radius $\delta(X)/4$.
Note that Theorem \ref{thm:grafting-differentials-close} bounds the
$L^1(X)$ norm, so the same upper bound applies to the
$L^1(D_{\delta/4})$ norm used here.

Since $(2 \phi_T(\lambda) + \phi_F(\lambda))$ is holomorphic, we apply
Lemma \ref{lem:extension} to the norm bound \eqref{eqn:partialbound},
and obtain
\begin{equation*}
\| 2 \phi_T(\lambda) + \phi_F(\lambda) \|_{L^1(X)} \leq C_3(X) \left (
1 + E(\lambda,X)^\half \right),
\end{equation*}
where $C_3(X)$ incorporates both the constants $C_1$ and $C_2$ from
\eqref{eqn:partialbound} and $C(\delta(X)/4,X)$ from Lemma
\ref{lem:extension}.  Since $\|\phi_F(\lambda)\|_{L^1(X)} =
E(\lambda,X)$, this proves Theorem \ref{thm:schwarzian-nearly-strebel}
for $\lambda$ supported on closed geodesics.

Finally, since both $\phi_T$ and $\phi_F$ are homeomorphisms, the
function $\lambda \mapsto \| 2 \phi_T(\lambda) + \phi_F(\lambda) \|_1$
is continuous on $\ML(S)$.  Since we have established a bound for this
function on the dense subset of $\ML(S)$ consisting of weighted simple
closed geodesics, the same bound applies to all laminations.  
\qedlabel{Theorem \ref{thm:schwarzian-nearly-strebel}}
\end{proof}

Note that the harmonic maps estimate (Theorem
\ref{thm:grafting-differentials-close}) provides the dominant factor in
the upper bound of Theorem \ref{thm:schwarzian-nearly-strebel}, and
the proof shows that any improvement to this estimate would give a
corresponding improvement to the bound on $\| 2 \phi_T(\lambda) +
\phi_F(\lambda) \|_1$:

\begin{cor}[of proof]
Fix $X \in \T(S)$.  For all $\lambda \in \ML(S)$ we have
\begin{equation*}
\| 2 \phi_T(\lambda) + \phi_F(\lambda) \|_{L^1(X)} \leq C(X) \left ( 1
+ \| \Phi(\lambda) - \phi_F(\lambda) \|_{L^1(X)} \right )
\end{equation*}
\end{cor}

It seems natural to ask whether it is possible to make the bound
completely independent of $\lambda$, that is:

\begin{question}
Is there a constant $C(X)$ such that $\| 2 \phi_T(\lambda) +
\phi_F(\lambda) \|_{L^1(X)} < C(X)$ for all $\lambda \in \ML(S)$?
\end{question}

Even if the bound can be made independent of $\lambda$, the dependence
on $X$ is probably necessary.  If there is a bound that depends only
on the topological type, one would need to select the right norm on
$Q(X)$ (perhaps using something other than $L^1$).  While all of the
standard norms on $Q(X)$ (e.g. $L^1$, $L^p$, $L^\infty$, etc.) are
equivalent when $X$ is fixed, the constants are not uniform as $X \to
\infty$.

\section{Applications: Holonomy and Fuchsian $\CP^1$ structures}
\label{sec:applications-fuchsian}

We now turn to an application of Theorem
\ref{thm:schwarzian-nearly-strebel} in the study of $\CP^1$ structures
and their holonomy representations.  Some background on this topic is
necessary before we state the results.

Recall that a projective surface $Z \in \P(S)$ has a holonomy
representation $\eta(Z) : \pi_1(S) \to \PSL_2(\C)$, which is unique up
to conjugation (see \S\ref{sec:grafting-and-cp1-structures}).  The
association of a holonomy representation to a projective structure
defines a map
\begin{equation*}
\eta : \P(S) \rightarrow \V(S)
\end{equation*}
where $\V(S)$ is the $\PSL_2(\C)$ character variety of
$\pi_1(S)$, i.e.
\begin{equation*}
\V(S) = \Hom(\pi_1(S), \PSL_2(\C)) /\!\!/ \PSL_2(\C)
\end{equation*}
and $\PSL_2(\C)$ acts on the set of homomorphisms by conjugation.
While this quotient (in the sense of geometric invariant theory) is
only defined up to birational equivalence, in this case one can
construct a good representative embedded in $\C^n$ using trace
functions (see \cite{culler-shalen}).  Then $\eta$ maps to the smooth
points of this variety and is a holomorphic local homeomorphism,
though it is not proper (and in particular is not a covering map).
The range of $\eta$ has been studied by a number of authors; recently,
it was shown that $\eta$ is essentially surjective onto a connected
component of $\V(S)$ \cite{GKM:monodromy}, proving a conjecture of
Gunning \cite{Gunning:projective-structures}.

We will be interested in the holonomy of $\CP^1$ structures on a fixed
Riemann surface $X$, i.e. the restriction of $\eta$ to $P(X)$.  The
resulting holomorphic immersion
\begin{equation*}
\eta_X : P(X) \to \V(S)
\end{equation*}
is proper \cite[\S 11.4]{GKM:monodromy} \cite{Tanigawa:divergence} and
injective \cite{Kra:generalization}.

Within $\V(S)$ there is the closed set $\AH(S)$ of discrete and
faithful representations, and its interior $\QF(S)$, which consists
of quasi-Fuchsian representations.  The quasi-Fuchsian representations
are exactly those whose limit sets are quasi-circles in $\CP^1$;
similarly, the set $\F(S)$ of Fuchsian representations consists of
those whose limit sets are round circles.

Let $K(X)$ denote the set of projective structures on $X$ with
discrete holonomy representations.  Shiga and Tanigawa showed that the
interior $\interior K(X)$ is exactly the set of projective
structures on $X$ with quasi-Fuchsian holonomy, i.e.
\begin{equation*}
\interior K(X) = \eta_X^{-1}(\QF(S)).
\end{equation*}
 
There is a natural decomposition of $\interior K(X)$ into countably
many open and closed subsets according to the topology of the
corresponding developing maps; the components in this decomposition
are naturally indexed by integral measured laminations (or
\emph{multicurves}) $\gamma \in \ML_\Z(S)$:
\begin{equation*}
\interior K(X) = \bigsqcup_{\gamma \in \ML_\Z(S)} B_\gamma(X)
\end{equation*}
The developing maps of projective structures in $B_\gamma(X)$ have
``wrapping'' behavior described by the lamination $\gamma$ (see
\cite[Ch. 7]{Kapovich:hyperbolic-manifolds}
\cite{Goldman:fuchsian-holonomy}).

The set $B_0(X)$ corresponding to the empty multicurve $0 \in
\ML_\Z(S)$ is the \emph{Bers slice} of $X$, i.e.~the set of projective
structures on $X$ that arise from the domains of discontinuity of
quasi-Fuchsian representations of $\pi_1(S)$
\cite{Shiga:projective-structures}.  When identified with a subset of
$Q(X)$ using the Schwarzian derivative, $B_0(X)$ is the connected
component of $\interior K(X)$ containing the origin, and is a bounded,
contractible open set.  The set $B_0(X)$ can also be characterized as
the set of projective structures on $X$ with quasi-Fuchsian holonomy
and injective developing maps.

We call the other sets $B_\gamma(X)$ (with $\gamma \neq 0$) the
\emph{exotic Bers slices}, because they are natural analogues of the
Bers slice $B_0(X)$, but they consist of quasi-Fuchsian projective
structures whose developing maps are not injective (which are called
``exotic'' $\CP^1$-structures).  Compared to $B_0(X)$, little is known
about the exotic Bers slices $B_\gamma(X)$; in particular it is not
known whether they are connected or bounded.

Each exotic Bers slice $B_\gamma(X)$ contains a distinguished point
$c_\gamma = c_\gamma(X)$, the \emph{Fuchsian center}, which is the
unique $\CP^1$ structure in $B_\gamma(X)$ with Fuchsian holonomy.
While in general the connection between the multicurve $\gamma$ and
the grafting laminations of projective structures in $B_\gamma(X)$ is
difficult to determine, for $c_\gamma$ the grafting lamination is just
$2 \pi \gamma$ (see \cite{Goldman:fuchsian-holonomy}), i.e.
\begin{equation*}
c_\gamma = \phi_T(2 \pi \gamma) \in B_\gamma(X).
\end{equation*}

The application of Theorem \ref{thm:schwarzian-nearly-strebel} we have
in mind involves the distribution of the Fuchsian centers within $P(X)
\simeq Q(X)$.  Associated to a multicurve $\gamma$ there is a
\emph{Jenkins-Strebel differential} $s_\gamma$,
\begin{equation*}
s_\gamma = \phi_F(2 \pi \gamma) \in Q(X)
\end{equation*}
which is a holomorphic quadratic differential whose noncritical
horizontal trajectories are closed and homotopic to the curves in
the support of $\gamma$.  Note that $-s_\gamma$ is then a differential
with closed vertical trajectories, so the sign in the definition is a
matter of convention.

Because the grafting laminations of the Fuchsian centers are integral
measured laminations (up to the factor $2 \pi$), Theorem
\ref{thm:schwarzian-nearly-strebel} implies that the associated points
in $Q(X)$ are close to differentials with closed trajectories.
Specifically, we have:

\begin{thm}
\label{thm:fuchsian-centers}
For any compact Riemann surface $X$ (with underlying smooth surface
$S$) and all multicurves $\gamma \in \ML_\Z(S)$, the Fuchsian center
$c_\gamma \in P(X) \simeq Q(X)$ and the $2 \pi$-integral
Jenkins-Strebel differential $s_\gamma \in Q(X)$ satisfy
\begin{equation*}
\| 2 c_\gamma + s_\gamma \|_1 \leq C(X) \left ( 1 +
\sqrt{\|s_\gamma\|_1} \right )
\end{equation*}
where $C(X)$ is a constant depending only on $X$.
\end{thm}

\begin{proof}
Using the definition of $c_\gamma$ and $s_\gamma$, this is immediate
from Theorem \ref{thm:schwarzian-nearly-strebel}:
\begin{equation*}
\| 2c_\gamma + s_\gamma \|_1 = \| 2\phi_T(2 \pi \gamma) + \phi_F(2 \pi
\gamma) \|_1 \leq C(X) \left ( 1 + \sqrt{\|\phi_F(2 \pi \gamma)
\|_{1}} \right )
\end{equation*}
\qedlabel{Theorem \ref{thm:fuchsian-centers}}
\end{proof}

Theorem \ref{thm:fuchsian-centers} implies that there is a rich
structure to the exotic Bers slices $B_\gamma(X)$ within the vector
space $Q(X)$ of holomorphic quadratic differentials; for example, each
rational ray $\R^+ \cdot s_\gamma \subset Q(X)$, $\gamma \in
\ML_Z(S)$ approximates the positions of an infinite sequence of
Fuchsian centers and their surrounding ``islands'' of quasi-Fuchsian
holonomy (though the distance from this line to the centers
may itself go to infinity, but at a slower rate).

For example, if we look at a sequence $\{ n \gamma \: | \: n =
1,2,\ldots\} \subset \ML_Z(S)$, then the norm $\| s_{n \gamma} \|_1 =
4 \pi^2 n^2 E(\gamma,X)$ grows quadratically with $n$, while by Theorem
\ref{thm:fuchsian-centers} we have $\| 2 c_{n \gamma} + s_{n
  \gamma} \|_1 = O(n)$ as $n \to \infty$.

\boldpoint{Numerical examples.}  In some cases the positions of the
Fuchsian centers in $Q(X)$ can be computed numerically.  While such
experiments do not yield new theoretical results, we present them here
to illustrate the connection between Fuchsian centers and
Jenkins-Strebel differentials (as in Theorem
\ref{thm:fuchsian-centers}) and the associated geometry of the
holonomy map. 

While the preceding discussion involved only compact surfaces, the
analogous theory of punctured surfaces with \emph{bounded} $\CP^1$
structures (those which are represented by meromorphic quadratic
differentials having at most simple poles at the punctures) is more
amenable to computation.  We will focus on the case where $X$ is a
punctured torus, so $P(X)$ is a one-dimensional complex affine space.
Here $X$ is commensurable with a planar Riemann surface (a
four-times-punctured $\CP^1$), so it is possible to compute the
holonomy representation of a $\CP^1$ structure on $X$ by numerical
integration of an ODE around contours in $\C$.  Existing discreteness
algorithms for punctured torus groups can then be used to create
pictures of Bers slices, the discreteness locus $K(X)$, and of the
positions of the Fuchsian centers within islands of quasi-Fuchsian
holonomy.  The images that follow were created using a computer
software package implementing these techniques \cite{Dumas:bear}. See
\cite{KS:bers-embedding-punctured-torus}
\cite{KSWY:drawing-bers-embedding} for further discussion of numerical
methods.

\begin{figure}
\fbox{\includegraphics[width=11cm]{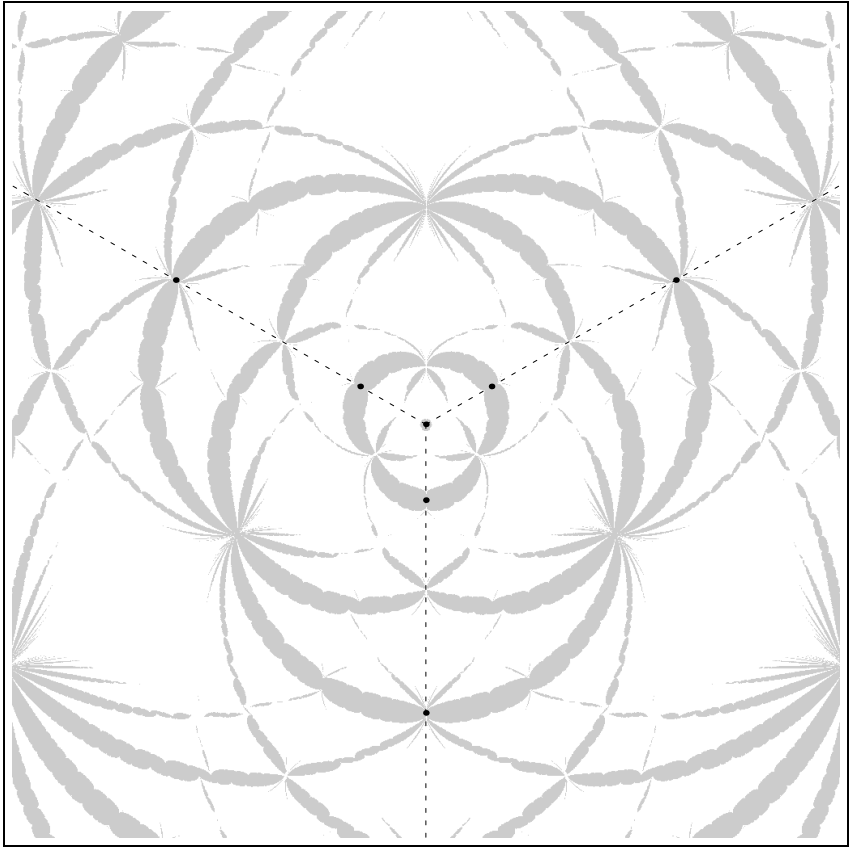}}
\caption{Fuchsian centers for the hexagonal punctured torus,
  corresponding to the empty multicurve (the center) and the first two
  multiples of the three systoles.}
\label{fig:fuchsian-centers-hex}
\end{figure}

Figure \ref{fig:fuchsian-centers-hex} shows part of $P(X)$ where $X$
is the hexagonal punctured torus, i.e. the result of identifying
opposite edges of a regular hexagon in $\C$ and removing one cycle of
vertices.  The image is centered on the Bers slice, and regions
corresponding to projective structures with discrete holonomy have
been shaded.  Each of the seven Fuchsian centers in this region of
$P(X)$ is marked; these correspond to the empty multicurve and the
three shortest curves on the hexagonal torus with multiplicities one
and two.

The dashed curves in Figure \ref{fig:fuchsian-centers-hex} are the
\emph{pleating rays} in $P(X)$ consisting of projective structures with
grafting laminations $\{ t \gamma \: | \: t \in \R^+ \}$, where
$\gamma$ is one of the three systoles.  Our use of the term ``pleating
ray'' is somewhat different than that of Keen-Series and others, as we
mean that the projective class of a bending lamination is fixed for a
family of equivariant pleated planes in $\H^3$, whereas in
\cite{KS:pleating-invariants} and elsewhere it is often assumed that
the associated $\PSL_2(\C)$ representation is quasi-Fuchsian and that
the pleated surface is one of its convex hull boundary surfaces.  The
pleating rays in $P(X)$ (as we have defined them) naturally
interpolate between the Fuchsian centers, which appear at the
$2\pi$-integral points.

The hexagonal torus is a special case because it has many symmetries.
In fact, the pleating rays for the systoles are forced (by symmetry)
to be Euclidean rays emanating from the origin in the directions of
the associated Jenkins-Strebel differentials.  Thus, while we expect
(based on Theorem \ref{thm:schwarzian-nearly-strebel} for compact
surfaces) that the Fuchsian centers associated to multiples of the
systoles lie near the lines of Jenkins-Strebel differentials, in this
example the Fuchsian centers lie on the lines, which coincide with the
pleating rays.

\begin{figure}
\fbox{\includegraphics[width=11cm]{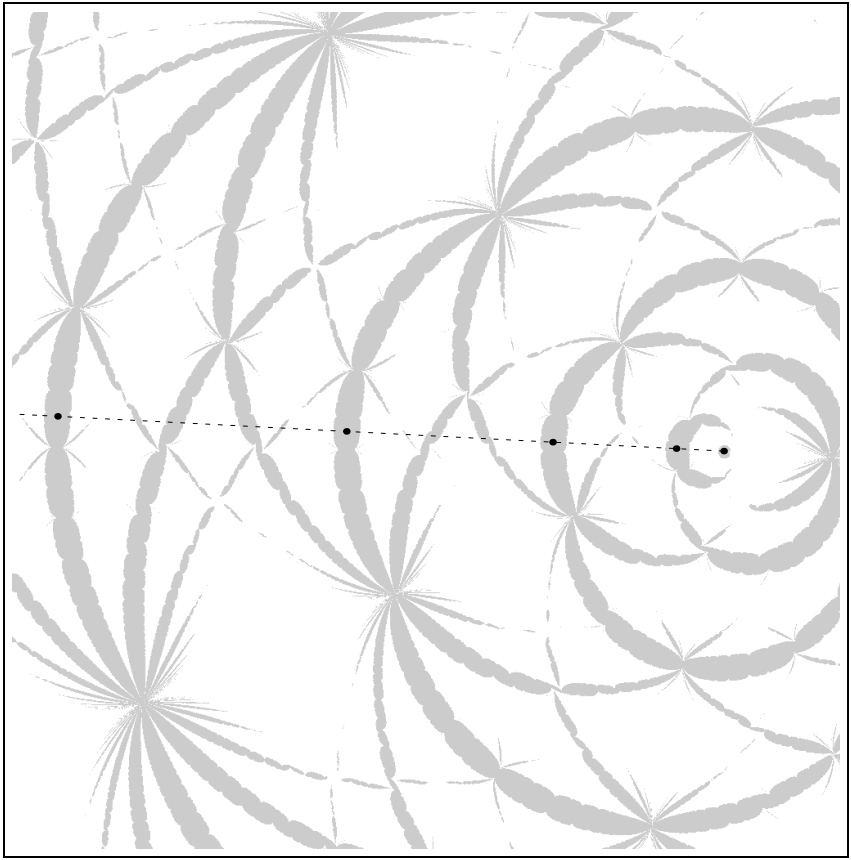}}
\caption{Fuchsian centers for a punctured torus with no symmetries;
  the modular parameter of the underlying Euclidean torus is
  approximately $\tau = 0.369 + 1.573 i$.  The right-most marked point
  is the origin, which is the Fuchsian center in the Bers slice.  The
  other Fuchsian centers correspond to multiples of a simple closed
  curve $\gamma$, which lie on the $\gamma$-pleating ray $\{ \phi_T(t
  \gamma) \: | \: t \in \R^+\}$ (the dashed curve).  The $L^1$
  distance between opposite sides of the image is approximately $250$.}
\label{fig:fuchsian-centers-nosym}
\end{figure}

Figure \ref{fig:fuchsian-centers-nosym} shows part of $P(X)$ for a
punctured torus $X$ without symmetries.  Five Fuchsian centers are
marked, corresponding to a certain simple closed curve $\gamma$ with
multiplicity $n$, $0 \leq n \leq 4$.  While there is no longer an
intrinsic symmetry that forces the centers to lie on a straight line
in $P(X) \simeq \C$, the pleating ray that contains these centers (the
dashed curve) is almost indistinguishable from the line of
Jenkins-Strebel differentials at the resolution of the figure.  In
fact, a sector of angle approximately $4 \times 10^{-5}$ centered on the
Jenkins-Strebel ray appears to contain the part of the pleating ray
visible in the figure.

\begin{figure}
\centerline{\includegraphics{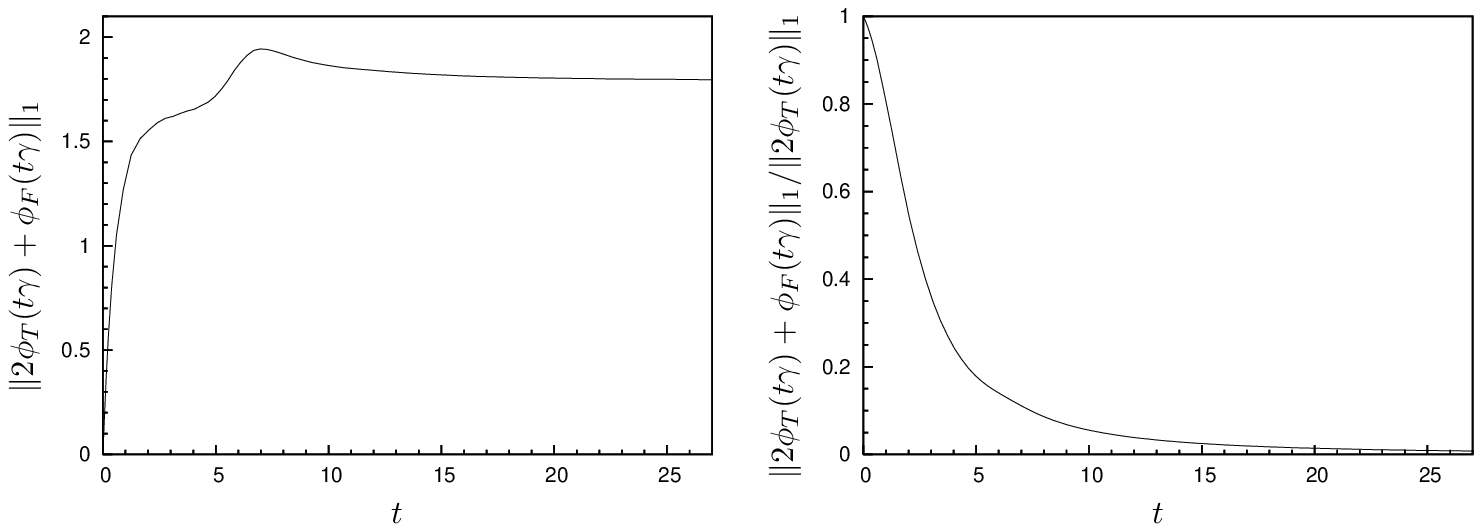}}
\caption{The $L^1$ norm of $(2 \phi_T(t \gamma) + \phi_F(t \gamma))
  \in Q(X)$ as a function of $t$ (left), and as a fraction of the
  $L^1$ norm of $2 \phi_T(t \gamma)$ (right).  Here $X$ is the
  punctured torus with modular parameter $\tau = 0.369 + 1.573 i$ and
  $\gamma$ is the curve on $X$ corresponding to the pleating ray in
  Figure \ref{fig:fuchsian-centers-nosym}.}
\label{fig:norm-difference}
\end{figure}

Figure \ref{fig:norm-difference} shows the norm of $2 \phi_T + \phi_F$
along the pleating ray from Figure \ref{fig:fuchsian-centers-nosym}.
In this region, this norm is everywhere less than $2$, while the first
Fuchsian center $c_\gamma$ has norm greater than $14$, and the width
of the region in Figure \ref{fig:fuchsian-centers-nosym} is
approximately $250$.  This should be contrasted with the distance
between the pleating ray and the line of Jenkins-Strebel differentials
as unparameterized curves, which is much smaller still.  Nevertheless,
the numerical results suggest that the Fuchsian centers are not
exactly collinear, with ratios of neighboring centers in Figure
\ref{fig:fuchsian-centers-nosym} apparently having small but nonzero
imaginary parts.

\section{Application: Compactification of $P(X)$ and $\P(S)$}
\label{sec:applications-compactification}

In this section we describe another application of Theorem
\ref{thm:schwarzian-nearly-strebel} that extends the results of
\cite{Dumas:antipodal} on $\CP^1$ structures and the asymptotics of
grafting and pruning.

As described in \S\ref{sec:grafting-and-cp1-structures}, Thurston's
projective extension of grafting provides a homeomorphism between the
space $\P(S)$ of $\CP^1$ structures on the differentiable compact
surface $S$ and the product $\ML(S) \times \T(S)$:
\begin{equation*}
\ML(S) \times \T(S) \xrightarrow{\Gr} \P(S).
\end{equation*}
Thus the set $P(X) \subset \P(S)$ of $\CP^1$ structures on a fixed
Riemann surface $X$ corresponds, using grafting, to a set of pairs
$M_X = \{ (\lambda,Y) \in \ML(S) \times \T(S) \: | \: \gr_\lambda Y =
X\}$.

In \cite{Dumas:antipodal}, it is shown that the image of $P(X)$ in
$\ML(S) \times \T(S)$ is well-behaved with respect to natural
compactifications of the two factors--the Thurston compactification
$\B{\T(S)}$ and the projective compactification $\B{\ML(S)} = \ML(S)
\sqcup \PML(S)$.  The asymptotic behavior of $P(X)$ is described in
terms of the \emph{antipodal involution} $i_X : \ML(S) \to \ML(S)$
(and its projectivization $i_X : \PML(S) \to \PML(S)$), where
$i_X(\lambda) = \mu$ if and only if $\lambda$ and $\mu$ are
measure-equivalent to the vertical and horizontal measured foliations
of a holomorphic quadratic differential $\phi \in Q(X)$
\cite[\S4]{Dumas:antipodal}.  Specifically, we have:

\begin{thm}[{\cite[Thm 1.2]{Dumas:antipodal}}]
\label{thm:antipodal}
For all $X \in \T(S)$, the boundary of $M_X = \{ (\lambda,Y) \in
\ML(S) \times \T(S) \: | \: \gr_\lambda Y = X\}$ in $\B{\ML(S)} \times
\B{\T(S)}$ is the graph of the projectivized antipodal involution, i.e.
\begin{equation*}
\B{M_X} = M_X \sqcup \Gamma(i_X)
\end{equation*}
where
\begin{equation*}
\Gamma(i_X) = \{ ( [\lambda], [i_X(\lambda)] ) \: | \: \lambda \in
\ML(S) \} \subset \PML(S) \times \PML(S).
\end{equation*}
\end{thm}

A deficiency in this description of the image of $P(X)$ in $\ML(S)
\times \T(S)$ is that it does not relate the boundary to the
Poincar\'e parameterization of $P(X)$ using quadratic differentials.

On the other hand, there \emph{is} a natural map from the boundary of the
space of holomorphic quadratic differentials to $\PML(S) \times
\PML(S)$: Consider the map
\begin{equation*}
\F^\perp \times \F : Q(X) \rightarrow \MF(S) \times \MF(S)
\end{equation*}
which records the vertical and horizontal measured foliations of a
holomorphic quadratic differential.  Using the natural identification
between $\MF(S)$ and $\ML(S)$, we also have a map $\Lambda^\perp \times
\Lambda : Q(X) \to \ML(S) \times \ML(S)$.  Because it is
homogeneous, there is an induced map between projective spaces
\begin{equation*}
\Lambda^\perp \times \Lambda : \PQ(X) \to \PML(S) \times \PML(S),
\end{equation*}
where $\PQ(X) = (Q(X) - \{ 0 \})/\R^+$.  The image of this map is, by
definition, the graph of the antipodal involution $i_X : \PML(S) \to
\PML(S)$, and $\Lambda^\perp \times \Lambda$ intertwines the action of
$-1$ on $\PQ(X)$ and the involution exchanging factors of $\PML(S)
\times \PML(S)$.

Since $P(X)$ and $Q(X)$ are identified using the Schwarzian
derivative, we can use the projective compactification $\B{Q(X)} =
Q(X) \sqcup \PQ(X)$ to obtain the \emph{Schwarzian compactification}
$\B{P(X)}_S$; a sequence of projective structures converges in
$\B{P(X)}_S$ if their Schwarzian derivatives can be rescaled (by
positive real factors) so as to converge in $( Q(X) - \{ 0 \} )$.

Using Theorem \ref{thm:schwarzian-nearly-strebel} we can extend
Theorem \ref{thm:antipodal} and show that the inclusion $P(X)
\into \ML(S) \times \T(S)$ extends continuously to
$\B{P(X)}_S$ and has $\Lambda^\perp \times \Lambda$ as its boundary
values:

\begin{thm}
\label{thm:schwarzian-compactification}
For each $X \in \T(S)$, the inclusion $P(X) \into \ML(S)
\times \T(S)$ obtained by taking Thurston's grafting coordinates for
the projective structures on $X$ extends continuously to a map
\begin{equation*}
\B{P(X)}_S \to \B{\ML(S)} \times \B{\T(S)}
\end{equation*}
where $\B{P(X)}_S = P(X) \sqcup \PQ(X)$ is the projective
compactification using the Schwarzian derivative, $\B{\T(S)}$ is the
Thurston compactification, and $\B{\ML(S)} = \ML(S) \sqcup \PML(S)$.
When restricted to the boundary, this extension agrees with the map
\begin{equation*}
\Lambda^\perp \times \Lambda : \PQ(X) \to \PML(S) \times \PML(S) 
\end{equation*}
which sends a projective class of differentials to the projective
measured laminations associated to its vertical and horizontal
foliations.  In particular, the boundary of $P(X)$ in $\B{\ML(S)}
\times \B{\T(S)}$ is the graph of the antipodal involution $i_X :
\PML(S) \to \PML(S)$.
\end{thm}

\begin{proof}
The idea is that the results of \cite{Dumas:antipodal} imply a similar
extension statement for a compactification of $P(X)$ using harmonic
maps, while Theorem \ref{thm:schwarzian-nearly-strebel} shows that
this compactification is the same as the one obtained using the
Schwarzian derivative (up to exchanging the vertical and horizontal
foliations).

Consider a divergent sequence in $P(X)$, with associated grafting
laminations $\lambda_i$; let $Y_i = \pr_{\lambda_i} X$.  The proof of
Theorem \ref{thm:antipodal} in \cite{Dumas:antipodal} uses Wolf's
theory of harmonic maps between Riemann surfaces and from Riemann
surfaces to $\R$-trees (see \cite{Wolf:high-energy-degeneration},
\cite{Wolf:harmonic-maps-to-r-trees}) to show that if $h_i : X \to
Y_i$ is the harmonic map compatible with the markings, and $\Phi_i \in
Q(X)$ is its Hopf differential, then
\begin{equation}
\begin{split}
\lim_{i \to \infty} [\Lambda(\Phi_i)] &= [\lambda] \in \PML(S)\\
\lim_{i \to \infty} [\Lambda^\perp(\Phi_i)] &= \lim_{i \to \infty} Y_i
= [i_X(\lambda)] \in \PML(S).
\end{split}
\label{eqn:hopf-limits}
\end{equation}
To rephrase these results, define the \emph{harmonic maps
compactification}
\begin{equation*}
\B{P(X)}_h = P(X) \sqcup \PQ(X)
\end{equation*}
where a sequence of projective structures converges to $[\Phi] \in
\PQ(X)$ if the sequence $\Phi_i$ of Hopf differentials of the
associated harmonic maps converges projectively to $\Phi$ (compare
\cite{Wolf:teichmuller-theory-harmonic-maps}).  Then
\eqref{eqn:hopf-limits} says that the inclusion $P(X) \into \ML(S)
\times \T(S)$ using the grafting coordinates extends continuously to
the harmonic maps compactification,
\begin{equation*}
\B{P(X)}_h \into \B{\ML(S)} \times \B{\T(S)},
\end{equation*}
and that the boundary of this extension is
\begin{equation*}
\Lambda \times \Lambda^\perp : \PQ(X) \to \PML(S) \times \PML(S).
\end{equation*}

By Theorem \ref{thm:schwarzian-nearly-strebel}, the projective limit
of the Schwarzian derivatives $\phi_T(\lambda_i)$ of a divergent
sequence $X(\lambda_i)$ is the same as that of the sequence
$-\phi_F(\lambda_i)$, where the sign is significant since $\PQ(X)$ is
the set of rays, rather than lines, in $Q(X)$.  On the other hand, by
Theorem \ref{thm:grafting-differentials-close}, the grafting
differentials $\Phi(\lambda_i)$ and $\phi_F(\lambda_i)$ have the same
projective limit, which is also the projective limit of the Hopf
differentials of the harmonic maps $h_i : X \to Y_{i}$ by Theorem 9.1
of \cite{Dumas:antipodal} (see also \cite{Dumas:erratum}).

As a result, the harmonic maps compactification and the Schwarzian
compactification are asymptotically related by the projectivization of
the linear map $-1 : \PQ(X) \to \PQ(X)$, which interchanges vertical and
horizontal foliations.  In particular the boundary map of $\B{P(X)}_S$
exists and is given by $\Lambda^\perp \times \Lambda = (\Lambda \times
\Lambda^\perp) \circ (-1)$.  \qedlabel{Theorem
\ref{thm:schwarzian-compactification}}
\end{proof}

\newcommand{\removethis}[1]{}

\end{document}